\documentclass[preprint,12pt]{elsarticle}
\usepackage{amsmath}
\usepackage{amsfonts}
\usepackage{amssymb}
\usepackage{graphicx}
\usepackage{geometry}
\usepackage{indentfirst}
\usepackage{subcaption}
\usepackage{float}
\usepackage{xcolor}

\usepackage{soul}
\usepackage{xcolor} 
\usepackage{algorithm}
\usepackage{algpseudocode}

\biboptions{numbers,sort&compress}

\journal{Additive Manufacturing}
\makeatletter
\let\ps@pprintTitle\ps@plain
\makeatother

\begin{document}

\begin{frontmatter}

\title{Localized Hybrid Analytical–Finite-Element Thermal Modeling for Laser Powder Bed Fusion}

\author[uw]{Mohammad Amin Mohammadi Golafshani}
\author[uw]{Krishnan Suresh}

\cortext[cor1]{Corresponding author}
\ead{ksuresh@wisc.edu}

\affiliation[uw]{
    organization={Department of Mechanical Engineering,
    University of Wisconsin--Madison},
    city={Madison},
    state={Wisconsin},
    country={USA}
}

\begin{abstract}
Hybrid analytical--finite-element methods support scan-resolved, part-scale thermal simulation of laser powder bed fusion, but their analytical evaluation cost grows as scan-line contributions accumulate. This work develops a space--time localized hybrid formulation that restricts where and how long each contribution is evaluated analytically. Residual terms and boundary conditions derived from the heat equation allow the finite-element correction to account for the thermal effect of contributions removed from the analytical field.
The formulation is assessed using a single-line benchmark against an independent full finite-element solution and multiline comparisons with the conventional hybrid formulation. Localization limits the number of analytical contributions retained while maintaining close agreement with the conventional solution. A comparison with uncompensated temporal truncation demonstrates the importance of the residual correction in retaining the thermal effect of removed contributions. In the 500-line case, analytical and total simulation times are reduced by factors of approximately \(25.5\) and \(3.92\), respectively, with a relative nodal 2-norm difference of \(0.139\%\).
\end{abstract}


\end{frontmatter}

\section{Introduction}

\begin{figure}[t]
    \centering
    \includegraphics[width=0.65\linewidth]{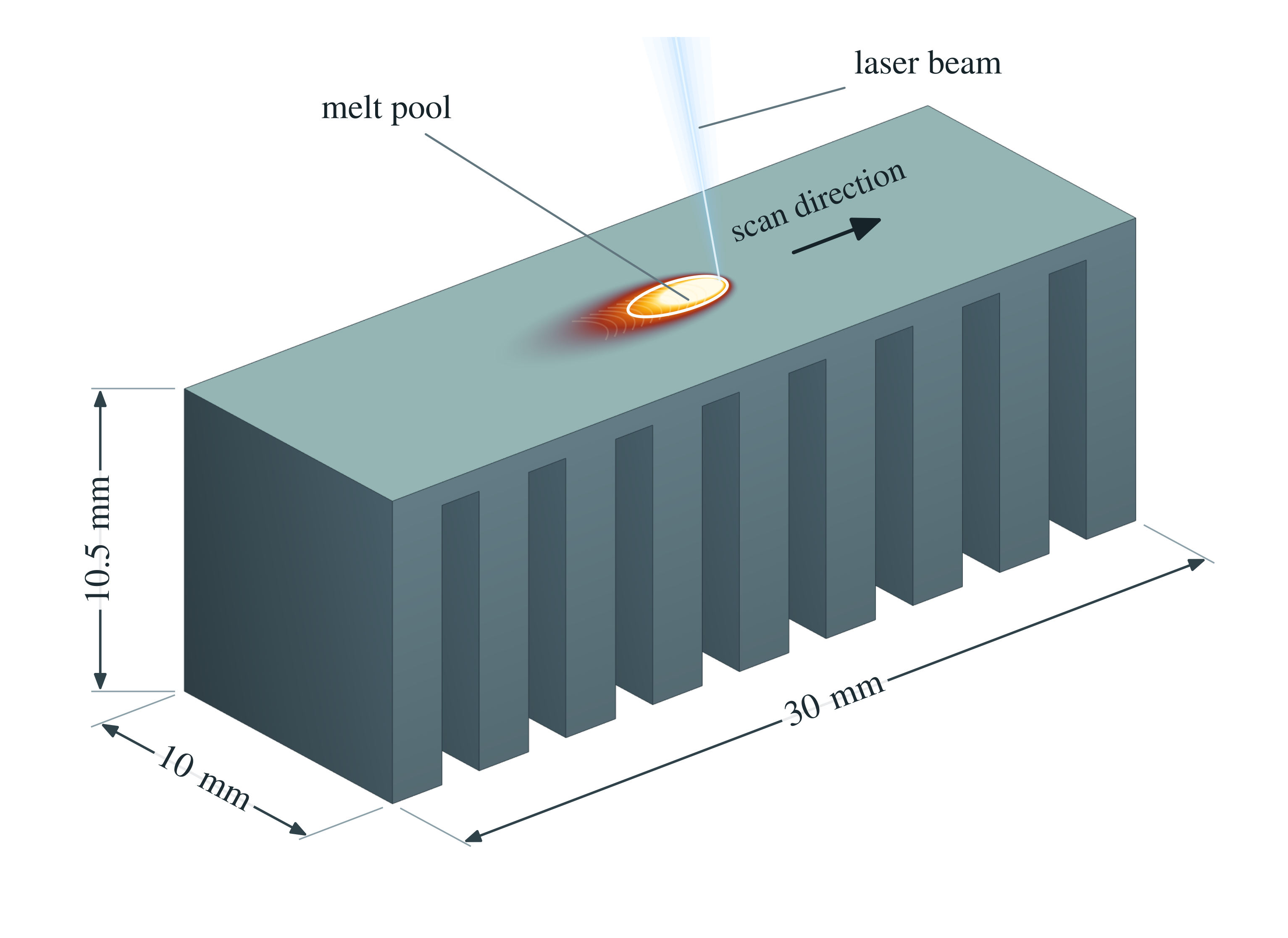}
    \caption{Representative fin geometry for part-scale LPBF thermal
    simulations, showing the moving laser and localized thermal response. The schematic was generated with assistance from
    Claude Opus 5 (Anthropic) and reviewed by the authors.}
    \label{fig:fin_partscale}
\end{figure}

Laser powder bed fusion (LPBF) is a metal additive manufacturing process
in which a high-energy laser selectively melts successive layers of metal
powder along a prescribed scan path (see Fig. \ref{fig:fin_partscale}). Its ability to fabricate geometrically
complex, high-accuracy components with high material efficiency has contributed
to its widespread adoption \cite{busachi2017review}. LPBF is increasingly used
to produce end-use components for aerospace, biomedical, automotive, and
energy applications \cite{kumar2016current,pereira2019comparison}.

LPBF involves complex thermal behavior arising from the highly localized
moving heat source and the associated rapid heating and cooling rates, which
produce steep temperature gradients within the material
\cite{promoppatum2017comprehensive, subedi2022towards, subedi2024deformation}. The resulting thermal history influences
melt-pool behavior, microstructural evolution, residual-stress formation, and
part distortion, thereby affecting the mechanical performance and structural
integrity of the fabricated component
\cite{soundararajan2021review,gong2014analysis}. Experimental characterization
of this thermal history is challenging because full-field temperature
measurements are limited, subsurface and rapidly evolving thermal fields are
difficult to capture, and extensive testing over large process-parameter
spaces is costly \cite{scheel2020adaptive,ji2022analytical}. Consequently,
predictive thermal simulations are essential for evaluating process conditions
and understanding part-scale thermal behavior beyond what can be measured
experimentally.

LPBF is a multiscale and multiphysics process in which different phenomena
dominate at different spatial and temporal scales, leading to distinct
modeling objectives and computational requirements
\cite{bresson2023global,de2025review,king2015overview,love2025advancing,
bayat2021review}. At the microscale, models examine microstructural phenomena
such as solidification, grain growth, and phase transformations, requiring
very fine spatial and temporal resolution
\cite{soundararajan2021review}. At the mesoscale, studies typically examine
melt-pool thermodynamics and fluid flow along individual scan tracks, together
with defect formation such as porosity evolution
\cite{xia2017porosity}. At the macroscopic or part scale, the objective shifts
to predicting the accumulated thermal history over the full build and its
thermomechanical consequences, including residual-stress accumulation and
part distortion
\cite{li2016multiscale,bayat2021review, prabhune2020fast}. To remain computationally feasible,
part-scale models employ simplified representations of localized melt-pool
physics. The present work therefore focuses on \emph{part-scale transient heat
conduction} rather than explicitly resolving fluid flow and phase-change
phenomena within the melt pool.

A wide range of approaches has been developed to predict thermal
histories in LPBF, with different trade-offs between physical
fidelity and computational efficiency. Fully numerical methods,
including finite-element analysis (FEA), can accommodate complex
geometries and nonlinear material behavior but may become
prohibitively expensive at the part scale because the localized
moving heat source requires fine spatial and temporal resolution
\cite{soundararajan2021review}. Analytical moving heat-source
solutions provide computationally efficient temperature predictions
but are generally restricted to idealized domains, simplified
material behavior, and a limited class of boundary conditions
\cite{eagar1983temperature}. Hybrid analytical--numerical methods combine these complementary capabilities by representing the heat-source response analytically and using a numerical correction to account for the finite part geometry and prescribed boundary conditions
\cite{yang2018semi,moran2021scan}. By capturing the steep
near-source temperature gradients analytically, these methods
can reduce the resolution required by the numerical mesh.
Different formulations use different analytical source solutions
and numerical methods to compute the correction.

The present work considers a hybrid formulation based on the
transient Eagar--Tsai solution~\cite{eagar1983temperature} and a finite-element correction.
After a scan line ends, the deposited heat continues to diffuse
and affect the temperature. As long as this effect is represented
by the line's analytical solution, that solution must be
reevaluated at later time steps. Retaining all earlier lines
therefore increases the analytical workload as scanning progresses.
Controlling this cost requires deciding where and for how long
each contribution is evaluated analytically, while accounting
for its thermal effect when it is removed. Related approaches to managing analytical contributions, including those developed within other formulations, are reviewed in Section~\ref{sec:litReview}.

To address this challenge, the present work introduces a
\emph{space--time localized analytical--finite-element formulation}.
Each analytical scan-line field is multiplied by smooth spatial
and temporal localization functions, restricting its explicit
evaluation to a neighborhood of the scan segment and a limited
time interval. Residual terms and boundary conditions derived
from the heat equation account for the changing analytical
weights through the finite-element correction. The resulting
formulation reduces analytical evaluation in space and time
while maintaining consistency with the governing heat equation.

The remainder of the paper is organized as follows. Section
\ref{sec:litReview} reviews numerical, analytical, and hybrid
approaches for LPBF thermal simulation and strategies for
managing analytical contributions. Section
\ref{sec:methodology} presents the analytical basis of the
hybrid method and develops the proposed space--time localized
formulation. Section \ref{sec:results} presents the numerical
validation and multiline studies. Finally, Section
\ref{sec:conclusions} summarizes the principal findings and
future directions.

\section{Literature Review}
\label{sec:litReview}

This section reviews the developments that motivate the present work. It first
examines part scale numerical thermal models and the strategies used to reduce
their computational cost. It then discusses analytical moving-source solutions and their
combination with numerical boundary corrections, including
hybrid analytical--finite-element formulations for finite
geometries. Finally, it considers the increasing cost of accumulated
analytical contributions and the approaches developed to manage the retained
scan history.

\subsection{Part-Scale Numerical Thermal Modeling}
\label{sec:part_scale_numerical_models}

Detailed LPBF models commonly use finite element, finite volume, or
computational fluid dynamics formulations and may include material
properties that vary with temperature, phase change, fluid flow in the
melt pool, and thermomechanical coupling
\cite{soundararajan2021review,de2025review}. Resolving the steep
thermal gradients around the rapidly moving laser requires a spatial
resolution comparable to the laser spot size and small time increments
along the scan path. Extending this resolution to a part containing
many tracks and layers results in very large computational requirements.
Detailed multiphysics simulations are therefore generally limited to
individual tracks, melt pools, or a small number of layers, while
simulations of complete parts rely on simplified physical descriptions
and specialized numerical strategies \cite{zhang2022review}.

Several reduction and acceleration strategies have been developed to
make thermal simulation of complete parts tractable. Adaptive remeshing
performed by layer maintains finer resolution near the currently
deposited material while using a coarser representation of previously
solidified regions \cite{olleak2020part}. Dugast et al. developed a GPU-based matrix-free finite-element solver that uses adaptive meshing in the build direction and excludes inactive elements from GPU memory transfer. They also examined layer lumping as a further cost-reduction strategy
\cite{dugast2021part}. Refinement can also be updated as scanning progresses rather than maintained across the entire active layer. Olleak et al. combined scanwise adaptive remeshing with a GPU-based matrix-free
finite-element solver, concentrating fine resolution around
selected scan tracks \cite{olleak2022enabling}. Liu et al. combined a characteristic time-based heat input model with adaptive mesh remapping to simulate thermal histories at track-scale resolution while reducing
computational cost \cite{liu2023macroscale}. These strategies reduce spatial or temporal discretization requirements or accelerate numerical calculations. For methods that use remeshing, these savings come with the additional work of managing mesh changes and transferring the temperature field between successive meshes.

Methods based on multiple resolution levels provide another way to
separate the fine thermal response around the laser from the smoother
temperature history across the part. In the GO-MELT framework, fine,
intermediate, and coarse discretizations are coupled and evaluated
using GPU acceleration \cite{leonor2024go}. Elahi et al. subsequently extended GO-MELT with temperature-dependent material properties and time-step subcycling between levels for part-scale LPBF thermal
simulations \cite{elahi2025efficient}. When the localized moving-source response is resolved directly, fine numerical resolution is still required near the laser, whether through mesh refinement or coupled discretizations.
An alternative is to represent this response analytically,
reducing the need for the numerical mesh to resolve the laser
length scale directly. Analytical and hybrid analytical--numerical
methods are reviewed next.

\subsection{Analytical and Hybrid Thermal Models}
\label{sec:analytical_hybrid_models}
Analytical thermal models provide an alternative way to predict
temperature fields generated by moving heat sources. Their development
originates from classical welding theory, particularly the moving point
and line source solutions introduced by Rosenthal
\cite{rosenthal1946theory}. Eagar and Tsai later considered a Gaussian
distribution of heat over the surface of a semi-infinite domain
\cite{eagar1983temperature}. This distribution accounts for the finite
width of the applied heat flux and provides a more suitable description
of a laser source than an ideal point source. Classical solutions
generally assume linear heat conduction, constant material properties,
and an infinite or semi-infinite domain. Their assumptions,
capabilities, and limitations in welding and related laser processes
are reviewed by Nasiri and Enzinger \cite{nasiri2019powerful}.

Analytical heat source models have subsequently been applied to LPBF
to predict temperature histories, melt pool dimensions, and cooling
rates under conduction-dominated conditions
\cite{zagade2021analytical,honarmandi2021rigorous}. Scan paths can be modeled by
superposing source contributions at successive positions or by
combining the contributions of moving sources along individual scan
lines. Wang et al. used analytical source solutions to model different
scan strategies and the temperature accumulated across multiple tracks
\cite{wang2023analytical}. Zagade et al. extended analytical modeling
to scan strategies with curved laser trajectories, including circular
and spiral paths \cite{zagade2024analytical}. Steuben et al. also
enriched a moving-source analytical solution to account approximately
for finite and complex part boundaries, temperature-dependent material
properties, and material accretion \cite{steuben2019enriched}.

A moving-source solution derived for a semi-infinite domain does
not, by itself, satisfy the boundary conditions of a finite part,
particularly when the laser approaches an external surface.
A common hybrid approach combines an analytical source field
with a numerical correction that enforces the prescribed
boundary conditions on the actual part. The analytical field
describes the response to the moving source in a simple
reference domain. By representing the steep temperature
gradients near the laser analytically, this approach allows
a coarser numerical mesh away from part boundaries, although
finer resolution may still be required near those boundaries.

Yang et al. introduced a semi-analytical formulation in which the
laser path was represented by a sequence of instantaneous point
sources~\cite{yang2018semi}. The associated analytical temperature
fields were combined with image fields and a smooth numerical
complementary field solved using finite differences to impose
the boundary conditions of the finite body. In later work,
the discrete point-source representation of each scan vector
was replaced by a continuous analytical line-source solution,
and the boundary conditions were enforced through a finite-element
correction with mesh refinement near surfaces
\cite{yang2020computationally}. This avoided the construction
and geometric restrictions of image sources.
Moran et al. developed a superposition-based finite-element
approach combining a transient Rosenthal-type solution for a
moving point heat source in a semi-infinite domain with a
finite-element correction to enforce the part boundary
conditions. The formulation was verified against a converged
conventional finite-element solution~\cite{moran2018utility}.
The method was subsequently applied to scan-by-scan simulation
of a realistic component while preserving the individual scan
vectors and their timing~\cite{moran2021scan}.

Recent developments have focused on improving the treatment of
boundaries and complex part geometries. Yang et al. introduced
modulated image sources to reduce the numerical resolution required
near curved adiabatic boundaries \cite{yang2025computational}.
More recently, they combined the analytical source representation
with an isogeometric discretization of the correction field,
allowing complex part boundaries to be represented accurately
with fewer numerical degrees of freedom
\cite{yang2026efficient}.

\subsection{Computational Cost and Management of Analytical Contributions}

Evaluating analytical heat-source fields can contribute
substantially to the computational cost of scan-resolved thermal
simulations. The workload depends on the source representation,
the spatial evaluation points, and the treatment of contributions
from earlier scans. If all scan contributions are retained and
reevaluated analytically, the work at each time step increases
as scanning progresses. Approaches to controlling this cost
include more efficient source evaluation, reuse of stored
solutions, and numerical representation of earlier thermal
contributions. These strategies are discussed below.

In the point-source formulation of Yang
et al.~\cite{yang2018semi}, evaluating the source and image fields
accounted for most of the reported computational time.
Sufficiently diffused source contributions and their associated
image fields were added to the numerical correction and removed
from subsequent source evaluation, using a transfer time based
on thermal diffusivity and grid spacing. The subsequent
continuous line-source formulation reduced the analytical cost
of evaluating each scan vector~\cite{yang2020computationally}.

Wolfer et al.~\cite{wolfer2019fast} used a different semi-analytical
formulation in which the source contribution over each time interval
was evaluated analytically and the previous temperature field was
propagated through Gaussian convolution. Stored line solutions were
reused, and analytical evaluation was restricted to the vicinity of
the laser. This avoided explicit superposition of the complete scan
history, but did not use a finite-element correction field.

Moran et al.~\cite{moran2021scan} reduced the retained analytical
history by projecting completed scan contributions onto the
finite-element correction once their temperature fields had
become sufficiently smooth to be represented by the coarse
mesh. These contributions were then removed from subsequent
analytical evaluation. The transfer was described as maintaining
thermal continuity, although the paper does not detail the
projection procedure or provide a quantitative criterion for
selecting the transfer time.

The present work develops a different treatment based on smooth
spatial and temporal localization functions. Rather than adding
or projecting completed analytical contributions onto a numerical
field at discrete transfer events, the formulation gradually
reduces their analytical weights and derives the corresponding
residual terms and boundary conditions from the heat equation.
Spatial localization also restricts where each retained
contribution is evaluated. The numerical studies examine the
resulting agreement and computational savings, while a comparison
with uncompensated temporal truncation demonstrates the
importance of accounting for removed contributions.

\section{Theory and Methodology}
\label{sec:methodology}

Consider a finite LPBF domain
\(\Omega\subset\mathbb{R}^{3}\) with boundary
\(\partial\Omega=\Gamma_D\cup\Gamma_N\), where
\(\Gamma_D\cap\Gamma_N=\varnothing\).
Under the assumptions of homogeneous and isotropic material properties,
constant thermophysical parameters, and heat transfer by conduction, the
temperature field \(T(\mathbf{x},t)\) satisfies
\begin{equation}
    \rho c_p \frac{\partial T}{\partial t}
    =
    k\nabla^2 T
    +
    Q(\mathbf{x},t),
    \qquad
    \mathbf{x}\in\Omega,
    \label{eq:governing_heat_equation}
\end{equation}
where \(\rho\), \(c_p\), and \(k\) denote the density, specific heat capacity, and thermal conductivity, respectively, \(Q(\mathbf{x},t)\) denotes the heat-source distribution associated with the moving laser, and \(\nabla^2\) is the Laplacian. For use in the hybrid formulation, define the linear heat operator
\begin{equation}
    \mathcal{L}(\,\cdot\,)
    :=
    \rho c_p
    \frac{\partial(\,\cdot\,)}{\partial t}
    -
    k\nabla^2(\,\cdot\,).
    \label{eq:heat_operator}
\end{equation}
Equation~\eqref{eq:governing_heat_equation} can then be written as
\begin{equation}
    \mathcal{L}T = Q.
    \label{eq:operator_governing}
\end{equation}

The governing equation is supplemented by the boundary and initial
conditions
\begin{align}
    T(\mathbf{x},t)
    &=
    T_b(\mathbf{x},t),
    && \mathbf{x}\in\Gamma_D,
    \label{eq:governing_dirichlet}
    \\
    -k\nabla T(\mathbf{x},t)\cdot\mathbf{n}
    &=
    \bar q(\mathbf{x},t),
    && \mathbf{x}\in\Gamma_N,
    \label{eq:governing_neumann}
    \\
    T(\mathbf{x},0)
    &=
    T_0(\mathbf{x}),
    && \mathbf{x}\in\Omega,
    \label{eq:governing_initial}
\end{align}
where \(\mathbf{n}\) is the outward unit normal,
\(T_b\) is the prescribed temperature on \(\Gamma_D\),
\(\bar q\) is the prescribed outward heat flux on \(\Gamma_N\), and
\(T_0\) is the initial temperature field.
In the present study, convective and radiative heat losses are neglected,
so the free Neumann boundaries are treated as insulated,
\(\bar q=0\). The laser heat input is represented through
\(Q(\mathbf{x},t)\) and is therefore not applied again through the
Neumann boundary condition.

The present study further neglects phase change and melt-pool fluid
flow since the focus is on part-scale thermal modeling. Together with the preceding assumptions, these simplifications
yield a linear thermal problem, allowing the temperature contributions
produced by individual laser exposures to be superposed.

The following subsections present the analytical basis and numerical
formulation of the proposed method. First, the Eagar--Tsai solution is
introduced for evaluating the temperature field generated by a moving
laser. Next, the localized analytical--finite-element formulation is
developed by introducing spatial and temporal localization into the
conventional hybrid decomposition. Finally, the weak form and its
finite-element discretization are derived for numerical implementation.

Throughout this section,
\(\mathbf{x}=(x,y,z)^{\mathsf{T}}\) denotes a spatial point. Boldface
symbols denote vectors and matrices, whereas scalar quantities are
written in standard font.

\subsection{Analytical solution for a moving Gaussian heat source}
\label{sec:analytical_solution}

The analytical component of the present formulation is based on the
classical Eagar--Tsai solution for a moving heat source in a homogeneous
semi-infinite solid \cite{eagar1983temperature}. The analytical domain
is bounded by a planar upper surface and extends infinitely in the depth
and lateral directions. The formulation assumes homogeneous and
isotropic material properties, constant thermophysical parameters, heat
transfer by conduction, and negligible convective and radiative heat
losses.

The laser is represented by a Gaussian heat-source distribution acting
on the upper surface. Let \(P_L\) denote the nominal laser power and
\(\eta\) the absorptivity. The total absorbed power is
\begin{equation}
    P_{\mathrm{abs}}
    =
    \eta P_L.
    \label{eq:absorbed_laser_power}
\end{equation}
For a source centered at the origin of the surface coordinate system,
the corresponding Gaussian surface power density is
\begin{equation}
    q_L(r)
    =
    \frac{P_{\mathrm{abs}}}{2\pi\sigma^2}
    \exp\!\left(
        -\frac{r^2}{2\sigma^2}
    \right),
    \label{eq:gaussian_source}
\end{equation}
where \(r\) is the in-plane radial distance from the source center and
\(\sigma\) is the standard deviation of the Gaussian distribution.
The distribution is normalized so that its integral over the upper
surface equals \(P_{\mathrm{abs}}\). This planar Gaussian source model
is used in the Eagar--Tsai solution and subsequent analytical LPBF
thermal models \cite{eagar1983temperature,wolfer2019fast,
wang2023analytical}.

Following Wolfer et al.~\cite{wolfer2019fast}, the source term in
Eq.~\eqref{eq:governing_heat_equation} can be written as

\begin{equation}
    Q(\mathbf{x},t)
    =
    q_L\!\left(
        \left\|
        \mathbf{x}_{\parallel}-\mathbf{x}_L(t)
        \right\|
    \right)
    \delta(z-z_{\mathrm{top}}),
    \label{eq:surface_supported_source}
\end{equation}
where \(\mathbf{x}_{\parallel}=(x,y)^{\mathsf{T}}\),
\(\mathbf{x}_L(t)\) is the in-plane laser position,
\(z_{\mathrm{top}}\) is the coordinate of the upper surface, and
\(\delta\) is the Dirac delta distribution. Equation
\eqref{eq:surface_supported_source} represents a planar Gaussian source
acting only on the upper surface rather than throughout the material
volume. 

Consider now a source that begins moving at time \(t_s\) from the
origin with constant velocity \(V\) in the positive \(x\)-direction.
The corresponding transient Eagar--Tsai solution gives the
temperature rise above the initial uniform temperature:
\begin{equation}
\begin{aligned}
T_{\mathrm{ET}}(x,y,z,t;t_s)
&=
\frac{P_{\mathrm{abs}}}
{\pi\rho c_p\sqrt{4\pi\alpha}}
\int_{t_s}^{t}
\frac{(t-t')^{-1/2}}
     {2\alpha(t-t')+\sigma^2}
\\
&\qquad\times
\exp\!\left(
-\frac{(x-V(t'-t_s))^2+y^2}
       {4\alpha(t-t')+2\sigma^2}
-\frac{z^2}
       {4\alpha(t-t')}
\right)\,dt',
\end{aligned}
\label{eq:eagar_tsai_solution}
\end{equation}
where
\begin{equation}
    \alpha
    =
    \frac{k}{\rho c_p}
    \label{eq:thermal_diffusivity}
\end{equation}
is the thermal diffusivity.

In the semi-infinite domain, the Eagar--Tsai temperature rise satisfies
the heat equation with the Gaussian surface source defined in
Eq.~\eqref{eq:surface_supported_source}:
\begin{equation}
    \mathcal{L}T_{\mathrm{ET}} = Q.
    \label{eq:eagar_tsai_operator}
\end{equation}
In this representation, the laser input is included in \(Q\), and the
Neumann condition on the planar surface is homogeneous.

Equation~\eqref{eq:eagar_tsai_solution} gives the temperature field
generated by a source that begins at \(t_s\) and remains active up to
the observation time \(t\). This expression is next used to construct
the contribution of a finite scan line with arbitrary position and
orientation.

Let scan line \(\ell\) start at time \(t_s\) from the point
\(\mathbf{S}\in\mathbb{R}^3\) and travel with constant speed \(V\) in
the unit direction \(\mathbf{e}\). Since the scan line lies on the
planar upper surface, let
\(\mathbf{R}\in\mathbb{R}^{3\times3}\) denote the in-plane rotation
that aligns the scan direction with the local \(\xi\)-axis, so that
\begin{equation}
    \mathbf{R}\mathbf{e}
    =
    (1,0,0)^{\mathsf{T}}.
\end{equation}
The local coordinates \((\xi,\eta,\zeta)\) of a global observation point
\(\mathbf{x}\) are defined by
\begin{equation}
    (\xi,\eta,\zeta)^{\mathsf{T}}
    :=
    \mathbf{R}(\mathbf{x}-\mathbf{S}).
    \label{eq:local_coords_S}
\end{equation}
In this coordinate system, the scan direction coincides with the positive
\(\xi\)-axis, while \(\zeta\) remains the depth coordinate. At time
\(t'\geq t_s\), the beam center therefore has the local coordinates
\begin{equation}
    \bigl(\xi_c(t'),\eta_c(t'),\zeta_c(t')\bigr)
    =
    \bigl(V(t'-t_s),0,0\bigr).
    \label{eq:beam_center_local}
\end{equation}

In the local coordinate system, the contribution of a source beginning
at \((\mathbf{S},t_s)\) and remaining active is denoted by \(T_1\) and
is obtained from Eq.~\eqref{eq:eagar_tsai_solution} as
\begin{equation}
    T_1(\mathbf{x},t)
    =
    T_{\mathrm{ET}}(\xi,\eta,\zeta,t;t_s).
    \label{eq:T1}
\end{equation}

For a finite scan line, the laser is active only over the interval
\begin{equation}
    t_s\leq t'\leq t_e,
    \qquad
    D:=t_e-t_s,
    \label{eq:exposure_window}
\end{equation}
where \(D\) is the exposure duration. The endpoint of the scan line is
\begin{equation}
    \mathbf{E}
    :=
    \mathbf{S}+VD\,\mathbf{e}.
    \label{eq:endpoint_E}
\end{equation}

Using the same rotation matrix \(\mathbf{R}\), define a second local
coordinate system with its origin at \(\mathbf{E}\):
\begin{equation}
    (\xi_E,\eta_E,\zeta_E)^{\mathsf{T}}
    :=
    \mathbf{R}(\mathbf{x}-\mathbf{E}).
    \label{eq:local_coords_E}
\end{equation}
Equations~\eqref{eq:local_coords_S} and
\eqref{eq:endpoint_E} give
\begin{equation}
    \xi_E=\xi-VD,
    \qquad
    \eta_E=\eta,
    \qquad
    \zeta_E=\zeta.
    \label{eq:coord_shift_relation}
\end{equation}

Now consider a fictitious source that starts at
\((\mathbf{E},t_e)\) and remains active for \(t'\geq t_e\). Its
temperature field is denoted by \(T_2\) and is given by
\begin{equation}
    T_2(\mathbf{x},t)
    =
    T_{\mathrm{ET}}(\xi_E,\eta_E,\zeta_E,t;t_e).
    \label{eq:T2}
\end{equation}

For \(t'\geq t_e\),
\[
    t'-t_e=(t'-t_s)-D.
\]
Using Eq.~\eqref{eq:coord_shift_relation},
\[
\begin{aligned}
\xi_E-V(t'-t_e)
&=
(\xi-VD)-V\bigl((t'-t_s)-D\bigr)
\\
&=
\xi-V(t'-t_s),
\end{aligned}
\]
with \(\eta_E=\eta\) and \(\zeta_E=\zeta\). Therefore, the kernels of
\(T_1\) and \(T_2\) are identical over the interval \([t_e,t]\), and
their difference removes the fictitious continuation of the laser after
the physical scan line has ended.

For \(t>t_e\), the temperature contribution of the finite scan line is
therefore
\begin{equation}
    T_\ell(\mathbf{x},t)
    =
    T_1(\mathbf{x},t)-T_2(\mathbf{x},t).
\end{equation}
Equivalently, a form that applies both while the line is being scanned
and after it has been completed is obtained by restricting the source
integration to the physical exposure interval. For \(t\geq t_s\),
\begin{equation}
\begin{aligned}
T_\ell(\mathbf{x},t)
&=
\frac{P_{\mathrm{abs}}}
{\pi\rho c_p\sqrt{4\pi\alpha}}
\int_{t_s}^{\,\min\{t,t_e\}}
\frac{(t-t')^{-1/2}}
     {2\alpha(t-t')+\sigma^2}
\\
&\qquad\times
\exp\!\left(
-\frac{(\xi-V(t'-t_s))^2+\eta^2}
       {4\alpha(t-t')+2\sigma^2}
-\frac{\zeta^2}
       {4\alpha(t-t')}
\right)\,dt'.
\end{aligned}
\label{eq:T_ell}
\end{equation}
Before the line starts,
\(T_\ell(\mathbf{x},t)=0\) for \(t<t_s\).
The contribution remains nonzero after the scan line is completed,
as the deposited heat continues to diffuse through the domain.

The quantity
\(T_\ell(\mathbf{x},t)\) in Eq.~\eqref{eq:T_ell} represents the
temperature contribution generated by scan line \(\ell\), and the
analytical temperature field associated with a complete scan path is
obtained by superposing the individual line contributions.

In addition to temperature, the spatial temperature gradient is required
in the hybrid formulation. Equation~\eqref{eq:T_ell} is evaluated in
the local coordinates \((\xi,\eta,\zeta)\), whereas the hybrid
formulation requires derivatives with respect to the global coordinates
\((x,y,z)\). From Eq.~\eqref{eq:local_coords_S}, the coordinate
Jacobian is \(\mathbf{R}\). Therefore, by the chain rule,
\begin{equation}
    \begin{pmatrix}
    \displaystyle \frac{\partial T_\ell}{\partial x}\\[2mm]
    \displaystyle \frac{\partial T_\ell}{\partial y}\\[2mm]
    \displaystyle \frac{\partial T_\ell}{\partial z}
    \end{pmatrix}
    =
    \mathbf{R}^{\mathsf{T}}
    \begin{pmatrix}
    \displaystyle \frac{\partial T_\ell}{\partial \xi}\\[2mm]
    \displaystyle \frac{\partial T_\ell}{\partial \eta}\\[2mm]
    \displaystyle \frac{\partial T_\ell}{\partial \zeta}
    \end{pmatrix}.
    \label{eq:gradient_transformation}
\end{equation}
Thus, the derivatives of Eq.~\eqref{eq:T_ell} are first evaluated in
the local coordinate system and then transformed by
\(\mathbf{R}^{\mathsf{T}}\) to obtain the temperature-gradient
components in the global coordinate system.

\subsection{Conventional hybrid analytical--finite-element formulation}
\label{sec:conventional_hybrid}

\begin{figure}[h]
    \centering
    \includegraphics[width=\textwidth]
    {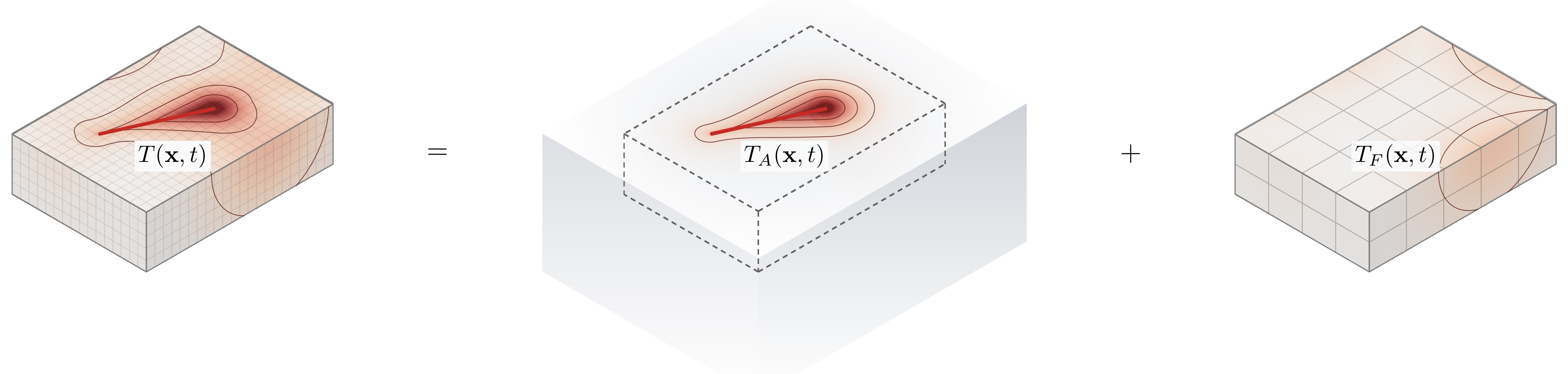}

    \makebox[\textwidth][c]{%
        \makebox[0.33\textwidth][c]{(a)}%
        \makebox[0.33\textwidth][c]{(b)}%
        \makebox[0.33\textwidth][c]{(c)}%
    }

    \caption{Schematic of the conventional hybrid decomposition.
    (a) Total temperature field in the finite part; the fine grid
    illustrates the spatial scale of the near-source response.
    (b) Analytical moving-source field in the semi-infinite reference
    domain.
    (c) Finite-element correction field computed on a coarser mesh
    in the finite part. The schematic was generated with assistance from
    Claude Opus 5 (Anthropic) and reviewed by the authors.}
    \label{fig:hybrid_method}
\end{figure}

Hybrid analytical--finite-element formulations combine the ability of
analytical moving-source solutions to capture steep temperature
gradients near the laser with the geometric flexibility of
finite-element analysis
\cite{moran2021scan,yang2020computationally,yang2025computational}.
The total temperature field is decomposed as
\begin{equation}
    T(\mathbf{x},t)
    =
    T_A(\mathbf{x},t)
    +
    T_F(\mathbf{x},t),
    \label{eq:conventional_hybrid}
\end{equation}
where \(T_A\) represents the laser-induced temperature field in a
semi-infinite reference domain and \(T_F\) is a finite-element
correction that enforces the boundary conditions of the finite part.
Figure~\ref{fig:hybrid_method} illustrates the decomposition and
the different spatial scales represented by the two fields.

For a scan path containing multiple lines, the analytical field is
obtained by superposition:
\begin{equation}
    T_A(\mathbf{x},t)
    =
    \sum_{\ell\in\mathcal{A}(t)}
    T_\ell(\mathbf{x},t),
    \label{eq:analytical_line_superposition}
\end{equation}
where \(\mathcal{A}(t)\) is the set of scan lines whose exposure has
begun by time \(t\). Let \(Q_\ell\) denote the laser-source contribution
associated with line \(\ell\). Each line field satisfies
\(\mathcal{L}T_\ell=Q_\ell\) in the semi-infinite reference domain.
Since the heat operator is linear, the analytical field satisfies
\begin{equation}
    \mathcal{L}T_A
    =
    \sum_{\ell\in\mathcal{A}(t)}Q_\ell
    =
    Q.
    \label{eq:analytical_operator}
\end{equation}

Substituting Eq.~\eqref{eq:conventional_hybrid} into the governing
equation \(\mathcal{L}T=Q\) and using
Eq.~\eqref{eq:analytical_operator} gives
\begin{equation}
    \mathcal{L}T_F
    =
    0.
    \label{eq:conventional_FE_equation}
\end{equation}

The correction field is solved in the finite domain \(\Omega\).
Its boundary and initial conditions follow from
Eq.~\eqref{eq:conventional_hybrid} and the prescribed conditions
for \(T\). On \(\Gamma_D\),
\begin{equation}
    T_F
    =
    T_b-T_A,
    \qquad \mathbf{x}\in\Gamma_D.
    \label{eq:conventional_dirichlet}
\end{equation}
On \(\Gamma_N\),
\begin{equation}
    -k\nabla T_F\cdot\mathbf{n}
    =
    \bar q+k\nabla T_A\cdot\mathbf{n},
    \qquad \mathbf{x}\in\Gamma_N.
    \label{eq:conventional_neumann}
\end{equation}
Initially,
\begin{equation}
    T_F(\mathbf{x},0)
    =
    T_0(\mathbf{x})-T_A(\mathbf{x},0),
    \qquad \mathbf{x}\in\Omega.
    \label{eq:conventional_initial}
\end{equation}
Thus, the analytical field represents the laser-source response,
while the finite-element field provides the correction required by
the boundary and initial conditions of the finite part.

Because the analytical field captures the steep temperature gradients
near the moving laser, the finite-element mesh needs to resolve only
the smoother correction field. This permits a coarser mesh than a
fully numerical model that directly resolves the laser-generated
field, reducing the number of finite-element unknowns. For explicit
time integration, a coarser mesh also relaxes the diffusion stability
restriction, which typically scales as
\(\Delta t\lesssim h^2/\alpha\), where \(h\) is a characteristic
element size \cite{strang2007computational}. For implicit schemes without this stability restriction,
the smoother correction can still permit larger time steps, provided
its evolution and the changing boundary conditions are adequately
resolved.

However, these computational advantages do not eliminate the
growing cost of analytical evaluation.
Equation~\eqref{eq:analytical_line_superposition} includes every
scan line whose exposure has begun. After a line is completed,
its analytical field continues to evolve as the deposited heat
diffuses. As long as this contribution remains in the analytical
representation, it must be reevaluated at subsequent time steps
wherever \(T_A\) or its gradient is required, including at points
far from the completed line. The analytical workload therefore
grows with the accumulated scan history. The localized
formulation developed next restricts these evaluations in space
and time and derives the correction terms needed to preserve
the thermal effect of the removed analytical contributions.

\subsection{Proposed Formulation}
\label{sec:localized_hybrid}

The proposed formulation restricts each analytical scan-line
contribution in space and time. This subsection first defines the
localization functions and then derives the corresponding correction
terms for the finite-element field.

\subsubsection{Spatial and temporal localization functions}

The conventional decomposition in
Eq.~\eqref{eq:conventional_hybrid} is generalized as
\begin{equation}
    T(\mathbf{x},t)
    =
    \sum_{\ell\in\mathcal{A}(t)}
    \chi_\ell(\mathbf{x},t)
    T_\ell(\mathbf{x},t)
    +
    T_F(\mathbf{x},t),
    \label{eq:localized_hybrid}
\end{equation}
where the localization function associated with scan line \(\ell\) is
written as
\begin{equation}
    \chi_\ell(\mathbf{x},t)
    =
    \chi_\ell^x(\mathbf{x})
    \chi_\ell^t(t),
    \qquad
    0\leq\chi_\ell^x\leq1,
    \qquad
    0\leq\chi_\ell^t\leq1.
    \label{eq:localization_function}
\end{equation}
Here, \(\chi_\ell^x(\mathbf{x})\) controls the spatial extent over which
the contribution of scan line \(\ell\) is explicitly represented,
whereas \(\chi_\ell^t(t)\) controls its temporal retention.

The localization functions are chosen to vary smoothly between their
fully active and inactive states. For the spatial localization associated
with scan line \(\ell\), let \(d_\ell(\mathbf{x})\) denote the shortest
distance from \(\mathbf{x}\) to the corresponding finite scan segment.
For a segment with start point \(\mathbf{x}_{s,\ell}\) and end point
\(\mathbf{x}_{e,\ell}\),
\begin{equation}
    d_\ell(\mathbf{x})
    =
    \min_{\lambda\in[0,1]}
    \left\|
        \mathbf{x}
        -
        \left[
            \mathbf{x}_{s,\ell}
            +
            \lambda
            \left(
                \mathbf{x}_{e,\ell}
                -
                \mathbf{x}_{s,\ell}
            \right)
        \right]
    \right\|.
    \label{eq:line_distance}
\end{equation}

The spatial localization function is defined as
\begin{equation}
    \chi_\ell^x(\mathbf{x})
    =
    \begin{cases}
        1,
        & d_\ell(\mathbf{x})\leq d_c,
        \\[1mm]
        S\!\left(
        \dfrac{d_\ell(\mathbf{x})-d_c}{w_x}
        \right),
        & d_c<d_\ell(\mathbf{x})<d_c+w_x,
        \\[3mm]
        0,
        & d_\ell(\mathbf{x})\geq d_c+w_x,
    \end{cases}
    \label{eq:spatial_localization}
\end{equation}
where \(d_c\) is the radius of the fully analytical core and \(w_x\)
is the width of the spatial transition region.

The temporal localization is defined relative to the completion time of
each scan-line contribution. Let \(t_{e,\ell}\) denote the completion
time of line \(\ell\), \(t_r\) the subsequent retention time, and
\(\Delta t_\chi\) the duration of the transfer to the finite-element
field. Then
\begin{equation}
    \chi_\ell^t(t)
    =
    \begin{cases}
        1,
        & t\leq t_{e,\ell}+t_r,
        \\[1mm]
        S\!\left(
        \dfrac{t-(t_{e,\ell}+t_r)}
              {\Delta t_\chi}
        \right),
        &
        t_{e,\ell}+t_r
        <
        t
        <
        t_{e,\ell}+t_r+\Delta t_\chi,
        \\[3mm]
        0,
        &
        t\geq
        t_{e,\ell}+t_r+\Delta t_\chi.
    \end{cases}
    \label{eq:temporal_localization}
\end{equation}

\begin{figure}[htbp]
    \centering

    \begin{subfigure}[t]{0.49\textwidth}
        \centering
        \includegraphics[width=\linewidth]{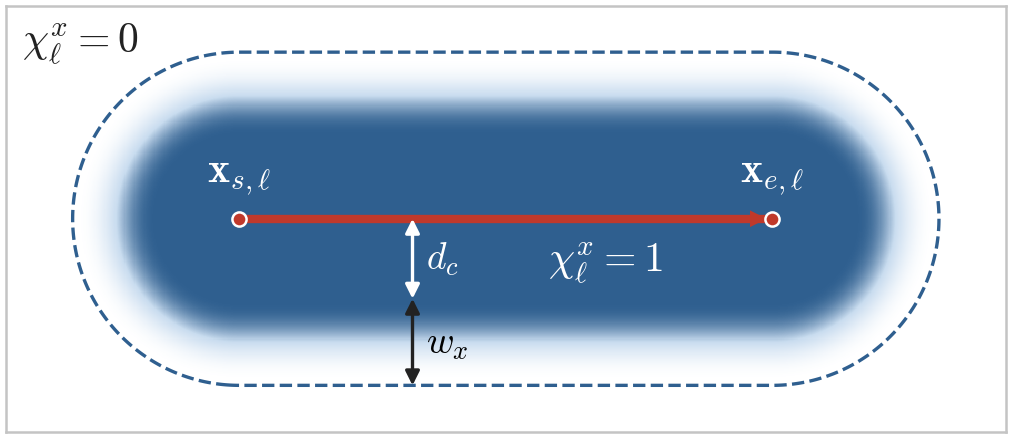}
        \caption{}
        \label{fig:spatial_localization_schematic}
    \end{subfigure}
    \hfill
    \begin{subfigure}[t]{0.49\textwidth}
        \centering
        \includegraphics[width=\linewidth]{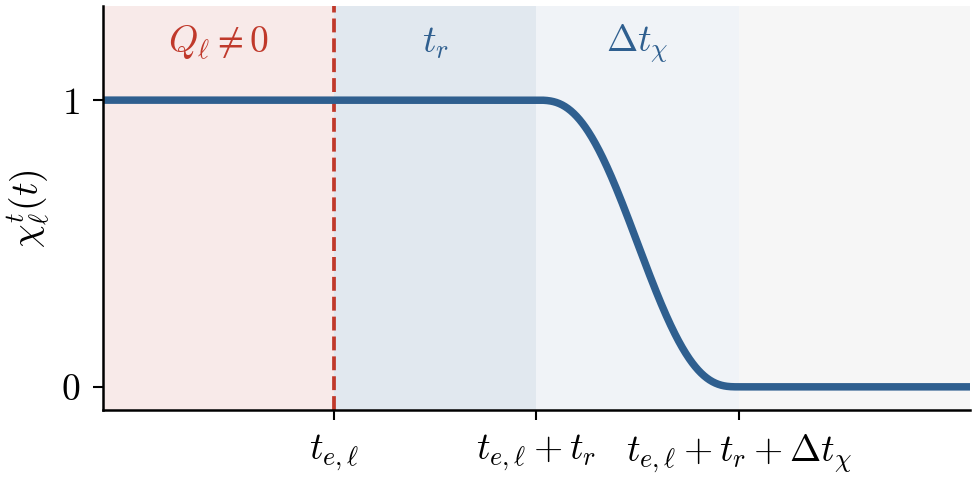}
        \caption{}
        \label{fig:temporal_localization_schematic}
    \end{subfigure}

    \caption{Spatial and temporal localization of an analytical
    scan-line contribution. (a) The spatial localization function
    satisfies \(\chi_\ell^x=1\) within the analytical core of radius
    \(d_c\), decreases smoothly over the transition width \(w_x\), and
    vanishes outside the localized region surrounding the scan segment.
    (b) The temporal localization function remains equal to one during
    the scan and for the retention interval \(t_r\) following the line
    completion time \(t_{e,\ell}\), and then decreases smoothly to zero
    over the transfer duration \(\Delta t_\chi\). The schematic was generated with assistance from Claude Opus 5 (Anthropic) and reviewed by the authors.}
    \label{fig:localization_functions}
\end{figure}

In Eqs.~\eqref{eq:spatial_localization} and
\eqref{eq:temporal_localization}, \(S:[0,1]\to[0,1]\) is a smooth,
monotonically decreasing transition function with \(S(0)=1\), \(S(1)=0\), and
\(S'(0)=S'(1)=0\). 
The spatial regions and temporal stages of localization are
shown schematically in Fig.~\ref{fig:localization_functions}.

The localized formulation is controlled by four parameters: the
analytical-core radius \(d_c\), the spatial transition width \(w_x\),
the post-completion retention time \(t_r\), and the temporal transfer
duration \(\Delta t_\chi\). During active scanning and throughout the
retention interval, \(\chi_\ell^t(t)=1\). Spatially, the analytical
contribution is fully retained for \(d_\ell\leq d_c\), smoothly reduced
for \(d_c<d_\ell<d_c+w_x\), and omitted from the explicit analytical
representation outside this region.
The conventional analytical field is recovered when
\(\chi_\ell^x=\chi_\ell^t=1\) for every scan-line contribution.

\subsubsection{Localization residual and finite-element correction}

To determine the governing equation for the correction field,
Eq.~\eqref{eq:localized_hybrid} is substituted into
Eq.~\eqref{eq:operator_governing}, giving
\begin{equation}
    \mathcal{L}T_F
    =
    Q
    -
    \sum_{\ell\in\mathcal{A}(t)}
    \mathcal{L}
    \left(
        \chi_\ell^x
        \chi_\ell^t
        T_\ell
    \right).
    \label{eq:localized_FE_preexpansion}
\end{equation}

Because \(\chi_\ell^x=\chi_\ell^x(\mathbf{x})\) and
\(\chi_\ell^t=\chi_\ell^t(t)\), the time derivative of one localized
analytical contribution is
\begin{equation}
    \frac{\partial
    \left(
        \chi_\ell^x\chi_\ell^tT_\ell
    \right)}
    {\partial t}
    =
    \chi_\ell^x
    \frac{d\chi_\ell^t}{dt}
    T_\ell
    +
    \chi_\ell^x
    \chi_\ell^t
    \frac{\partial T_\ell}{\partial t},
    \label{eq:localized_time_derivative}
\end{equation}
Furthermore,
\begin{equation}
\begin{aligned}
    \nabla^2
    \left(
        \chi_\ell^x
        \chi_\ell^t
        T_\ell
    \right)
    ={}&
    \chi_\ell^t
    \chi_\ell^x
    \nabla^2T_\ell
    \\
    &+
    2\chi_\ell^t
    \nabla\chi_\ell^x
    \cdot
    \nabla T_\ell
    +
    \chi_\ell^t
    T_\ell
    \nabla^2\chi_\ell^x .
\end{aligned}
\label{eq:localized_laplacian}
\end{equation}
Therefore,
\begin{align}
    \mathcal{L}
    \left(
        \chi_\ell^x
        \chi_\ell^t
        T_\ell
    \right)
    ={}&
    \chi_\ell^x
    \chi_\ell^t
    \left(
        \rho c_p
        \frac{\partial T_\ell}{\partial t}
        -
        k\nabla^2T_\ell
    \right)
    \nonumber\\
    &+
    \rho c_p
    \chi_\ell^x
    T_\ell
    \frac{d\chi_\ell^t}{dt}
    \nonumber\\
    &-
    2k\chi_\ell^t
    \nabla\chi_\ell^x
    \cdot
    \nabla T_\ell
    \nonumber\\
    &-
    k\chi_\ell^t
    T_\ell
    \nabla^2\chi_\ell^x.
    \label{eq:localized_operator_expansion}
\end{align}
Since \(\mathcal{L}T_\ell=Q_\ell\), this becomes
\begin{align}
    \mathcal{L}
    \left(
        \chi_\ell^x
        \chi_\ell^t
        T_\ell
    \right)
    ={}&
    \chi_\ell^x
    \chi_\ell^t
    Q_\ell
    +
    \rho c_p
    \chi_\ell^x
    T_\ell
    \frac{d\chi_\ell^t}{dt}
    \nonumber\\
    &-
    2k\chi_\ell^t
    \nabla\chi_\ell^x
    \cdot
    \nabla T_\ell
    \nonumber\\
    &-
    k\chi_\ell^t
    T_\ell
    \nabla^2\chi_\ell^x .
    \label{eq:localized_operator_final}
\end{align}

Substituting this expression into
Eq.~\eqref{eq:localized_FE_preexpansion} gives
\begin{align}
    \rho c_p
    \frac{\partial T_F}{\partial t}
    -
    k\nabla^2T_F
    ={}&
    \sum_{\ell\in\mathcal{A}(t)}
    \Big[
    \left(
        1-
        \chi_\ell^x\chi_\ell^t
    \right)Q_\ell
    \nonumber\\
    &\qquad
    -
    \rho c_p
    \chi_\ell^x
    T_\ell
    \frac{d\chi_\ell^t}{dt}
    \nonumber\\
    &\qquad
    +
    2k\chi_\ell^t
    \nabla\chi_\ell^x
    \cdot
    \nabla T_\ell
    \nonumber\\
    &\qquad
    +
    k\chi_\ell^t
    T_\ell
    \nabla^2\chi_\ell^x
    \Big].
    \label{eq:localized_FE_general}
\end{align}

When the laser is not scanning line \(\ell\), \(Q_\ell=0\),
so the source-localization term
\((1-\chi_\ell^x\chi_\ell^t)Q_\ell\) is zero. During exposure,
\(\chi_\ell^t=1\). The term is then also zero inside the
spatial core, where \(\chi_\ell^x=1\).
The only remaining case is active exposure at points where
\(\chi_\ell^x<1\). The beam center lies on the scan segment, so
these points are at least \(d_c\) from the beam center. From
Eq.~\eqref{eq:gaussian_source}, the laser power density there
is at most \(\exp[-d_c^2/(2\sigma^2)]\) times its value at
the beam center. For the parameter values used in this work,
\(\sigma=50~\mu\mathrm{m}\). The selected core radius
\(d_c=1~\mathrm{mm}\) gives a factor of
\(e^{-200}\approx1.38\times10^{-87}\). Even the smallest
core radius examined, \(d_c=0.5~\mathrm{mm}\), gives only
\(e^{-50}\approx1.93\times10^{-22}\). Thus, where the
source-localization term is not exactly zero, the Gaussian
source it contains is extremely small. Dropping this term gives
\begin{align}
    \rho c_p
    \frac{\partial T_F}{\partial t}
    -
    k\nabla^2T_F
    ={}&
    \sum_{\ell\in\mathcal{A}(t)}
    \Big[
    -
    \rho c_p
    \chi_\ell^x
    T_\ell
    \frac{d\chi_\ell^t}{dt}
    \nonumber\\
    &\qquad
    +
    2k\chi_\ell^t
    \nabla\chi_\ell^x
    \cdot
    \nabla T_\ell
    \nonumber\\
    &\qquad
    +
    k\chi_\ell^t
    T_\ell
    \nabla^2\chi_\ell^x
    \Big].
    \label{eq:localized_FE}
\end{align}

The right-hand side of Eq.~\eqref{eq:localized_FE} represents the
correction introduced by localization of the analytical field.
For each scan-line contribution, it can be separated into temporal and
spatial components,
\begin{equation}
    R_{\mathrm{loc},\ell}
    =
    R_{t,\ell}
    +
    R_{x,\ell},
    \label{eq:localization_residual}
\end{equation}
where
\begin{equation}
    R_{t,\ell}
    =
    -\rho c_p
    \chi_\ell^x
    T_\ell
    \frac{d\chi_\ell^t}{dt},
    \label{eq:temporal_residual}
\end{equation}
and
\begin{equation}
    R_{x,\ell}
    =
    2k\chi_\ell^t
    \nabla\chi_\ell^x
    \cdot
    \nabla T_\ell
    +
    k\chi_\ell^t
    T_\ell
    \nabla^2\chi_\ell^x.
    \label{eq:spatial_residual}
\end{equation}
The complete localization residual is
\begin{equation}
    R_{\mathrm{loc}}
    =
    \sum_{\ell\in\mathcal{A}(t)}
    R_{\mathrm{loc},\ell}.
\end{equation}
The temporal term \(R_{t,\ell}\) accounts for the gradual removal of an
analytical contribution from the explicit analytical representation.
The spatial term \(R_{x,\ell}\) accounts for the transition between
regions in which that analytical contribution is fully represented and
regions in which its effect is represented through the finite-element
correction field.

The boundary conditions for \(T_F\) follow from the decomposition in
Eq.~\eqref{eq:localized_hybrid}. On the Dirichlet boundary,
\begin{equation}
    T_F
    =
    T_b
    -
    \sum_{\ell\in\mathcal{A}(t)}
    \chi_\ell^x
    \chi_\ell^t
    T_\ell,
    \qquad
    \mathbf{x}\in\Gamma_D.
    \label{eq:localized_dirichlet}
\end{equation}
On the Neumann boundary,
\begin{equation}
\begin{aligned}
    -k\nabla T_F\cdot\mathbf{n}
    ={}&
    \bar q
    +
    k
    \sum_{\ell\in\mathcal{A}(t)}
    \chi_\ell^t
    \left(
        \chi_\ell^x\nabla T_\ell
        +
        T_\ell\nabla\chi_\ell^x
    \right)
    \cdot\mathbf{n},
    \\
    &\hspace{47mm}
    \mathbf{x}\in\Gamma_N.
\end{aligned}
\label{eq:localized_neumann}
\end{equation}

Together, the localization residual and the modified boundary conditions
allow the finite-element field to account for contributions removed from
the explicit analytical representation. The analytical evaluation can
therefore remain concentrated near recent scan lines, limiting its
spatial workload and the amount of scan history retained as scanning
progresses.

To obtain the weak form, define the test space
\[
    V_0
    :=
    \left\{
        v\in H^1(\Omega):
        v|_{\Gamma_D}=0
    \right\}.
\]
Multiplying Eq.~\eqref{eq:localized_FE} by \(v\in V_0\),
integrating over \(\Omega\), and integrating the diffusion
term by parts gives
\begin{equation}
\begin{aligned}
    \int_\Omega \rho c_p v
    \frac{\partial T_F}{\partial t}\,d\mathbf{x}
    +
    \int_\Omega k\nabla v\cdot\nabla T_F\,d\mathbf{x}
    ={}&
    \int_\Omega vR_{\mathrm{loc}}\,d\mathbf{x}
    \\
    &+
    \int_{\Gamma_N}kv\nabla T_F\cdot\mathbf{n}\,ds.
\end{aligned}
\label{eq:localized_weak_preliminary}
\end{equation}

The Laplacian of the spatial localization function in
Eq.~\eqref{eq:spatial_residual} is also integrated by parts.
Using \(\nabla(vT_\ell)=T_\ell\nabla v+v\nabla T_\ell\) gives
\begin{equation}
\begin{aligned}
    \int_\Omega vR_{x,\ell}\,d\mathbf{x}
    ={}&
    k\chi_\ell^t
    \int_\Omega
    \left(
        v\nabla\chi_\ell^x\cdot\nabla T_\ell
        -
        T_\ell\nabla v\cdot\nabla\chi_\ell^x
    \right)d\mathbf{x}
    \\
    &+
    k\chi_\ell^t
    \int_{\Gamma_N}
    vT_\ell\nabla\chi_\ell^x\cdot\mathbf{n}\,ds.
\end{aligned}
\label{eq:spatial_residual_weak}
\end{equation}
The boundary integral in Eq.~\eqref{eq:spatial_residual_weak}
cancels the term containing \(T_\ell\nabla\chi_\ell^x\) when
the Neumann condition in Eq.~\eqref{eq:localized_neumann}
is substituted into Eq.~\eqref{eq:localized_weak_preliminary}.
The resulting weak form, valid for every \(v\in V_0\), is
\begin{equation}
\begin{aligned}
    &\int_\Omega \rho c_p v
    \frac{\partial T_F}{\partial t}\,d\mathbf{x}
    +
    \int_\Omega k\nabla v\cdot\nabla T_F\,d\mathbf{x}
    \\
    &\quad=
    \sum_{\ell\in\mathcal{A}(t)}
    \int_\Omega
    \Bigg[
        -\rho c_p v\chi_\ell^xT_\ell
        \frac{d\chi_\ell^t}{dt}
        +
        k\chi_\ell^t
        \left(
            v\nabla\chi_\ell^x\cdot\nabla T_\ell
            -
            T_\ell\nabla v\cdot\nabla\chi_\ell^x
        \right)
    \Bigg]d\mathbf{x}
    \\
    &\qquad-
    \int_{\Gamma_N}v\bar q\,ds
    -
    k\sum_{\ell\in\mathcal{A}(t)}
    \int_{\Gamma_N}
    v\chi_\ell^x\chi_\ell^t
    \nabla T_\ell\cdot\mathbf{n}\,ds.
\end{aligned}
\label{eq:localized_weak_form}
\end{equation}
Thus, the weak form accounts for the spatial Laplacian through
first derivatives and does not require its explicit evaluation.

The correction field is approximated as
\(T_F(\mathbf{x},t)=\sum_jN_j(\mathbf{x})\theta_{F,j}(t)\),
where \(N_j\) are the finite-element shape functions.
Taking \(v=N_i\), with \(N_i=0\) on \(\Gamma_D\), gives
\begin{equation}
    \mathbf{C}\dot{\boldsymbol{\theta}}_F
    +
    \mathbf{K}\boldsymbol{\theta}_F
    =
    \mathbf{f}_{\mathrm{loc}}
    +
    \mathbf{f}_{\Gamma},
    \label{eq:localized_FE_discrete}
\end{equation}
where \(\boldsymbol{\theta}_F\) contains the nodal values of the
correction field, and
\begin{equation}
\begin{aligned}
    C_{ij}
    &=
    \int_\Omega \rho c_p N_iN_j\,d\mathbf{x},
    \\
    K_{ij}
    &=
    \int_\Omega k\nabla N_i\cdot\nabla N_j\,d\mathbf{x},
    \\
    (f_{\mathrm{loc}})_i
    &=
    \sum_{\ell\in\mathcal{A}(t)}
    \int_\Omega
    \Bigg[
        -\rho c_p N_i\chi_\ell^xT_\ell
        \frac{d\chi_\ell^t}{dt}
        +
        k\chi_\ell^t
        \left(
            N_i\nabla\chi_\ell^x\cdot\nabla T_\ell
            -
            T_\ell\nabla N_i\cdot\nabla\chi_\ell^x
        \right)
    \Bigg]d\mathbf{x},
    \\
    (f_{\Gamma})_i
    &=
    -\int_{\Gamma_N}N_i\bar q\,ds
    -
    k\sum_{\ell\in\mathcal{A}(t)}
    \int_{\Gamma_N}
    N_i\chi_\ell^x\chi_\ell^t
    \nabla T_\ell\cdot\mathbf{n}\,ds.
\end{aligned}
\label{eq:localized_FE_assembly}
\end{equation}
The Dirichlet condition is imposed using
Eq.~\eqref{eq:localized_dirichlet}.
Algorithm~\ref{alg:localized_hybrid} summarizes the resulting
time-stepping procedure.

\begin{algorithm}[H]
\caption{Localized hybrid thermal simulation}
\label{alg:localized_hybrid}
\begin{algorithmic}[1]
\Require Scan path, finite-element mesh, and localization parameters
\(d_c,w_x,t_r,\Delta t_\chi\)
\State Initialize \(T_F\) from the prescribed initial temperature.
\For{each time step \(t_n\)}
    \State Select lines that have started and satisfy
    \(\chi_\ell^t(t_n)>0\).
    \State Initialize the analytical field and load vectors.
    \For{each selected line \(\ell\)}
        \State Compute \(\chi_\ell^t\) and
        \(d\chi_\ell^t/dt\).
        \State Evaluate \(\chi_\ell^x\), \(\nabla\chi_\ell^x\),
        \(T_\ell\), and \(\nabla T_\ell\) where needed.
        \State Add \(\chi_\ell^x\chi_\ell^tT_\ell\)
        to the analytical field.
        \State Assemble the line's contributions to
        \(\mathbf{f}_{\mathrm{loc}}\) and \(\mathbf{f}_{\Gamma}\)
        using Eq.~\eqref{eq:localized_FE_assembly}.
    \EndFor
    \State Include the prescribed flux contribution in
    \(\mathbf{f}_{\Gamma}\) using
    Eq.~\eqref{eq:localized_FE_assembly}.
    \State Form the finite-element system in
    Eq.~\eqref{eq:localized_FE_discrete} and impose
    Eq.~\eqref{eq:localized_dirichlet}.
    \State Advance \(T_F\) and reconstruct \(T\) using
    Eq.~\eqref{eq:localized_hybrid}.
\EndFor
\end{algorithmic}
\end{algorithm}

ChatGPT 6 (OpenAI) was used to assist in developing the
Python simulation code. The implementation was tested
and validated by the authors.

\section{Results and discussion}
\label{sec:results}
The proposed localized hybrid formulation is evaluated on a fin
geometry (see Fig.~\ref{fig:fin_partscale}) through four studies.
First, the conventional and localized hybrid solutions are compared
against an independent full finite-element solution obtained using
Abaqus for a single-line benchmark. Next, a comparison with
uncompensated temporal truncation examines the importance of
accounting for removed analytical contributions. The effect of
localization parameters on accuracy and computational cost is
then examined. Finally, the analytical history and computational
savings are analyzed for a 500-line case. The last three studies
use the conventional hybrid formulation as the reference.

The fin geometry (see Fig. \ref{fig:fin_partscale}) is assigned homogeneous and isotropic material properties with constant mass density,
thermal conductivity, and specific heat capacity. The initial temperature is set to \(32~^\circ\mathrm{C}\) throughout
the domain, and the bottom surface is maintained at this temperature
during the simulation. The remaining external surfaces are treated
as adiabatic. The laser power,
absorptivity, beam size, and scan speed are also kept fixed throughout
the numerical studies. Their values are listed in
Table~\ref{tab:physical_parameters}.

\begin{table}[htbp]
    \centering
    \caption{Physical parameters common to all numerical studies.}
    \label{tab:physical_parameters}
    \begin{tabular}{lll}
        \hline
        Category & Parameter & Value \\
        \hline
        Material
        & Mass density \(\rho\)
        & \(8440~\mathrm{kg\,m^{-3}}\) \\

        & Thermal conductivity \(k\)
        & \(11.4~\mathrm{W\,m^{-1}\,K^{-1}}\) \\

        & Specific heat capacity \(c_p\)
        & \(443~\mathrm{J\,kg^{-1}\,K^{-1}}\) \\
        \hline
        Laser and scan
        & Laser power \(P_L\)
        & \(300~\mathrm{W}\) \\

        & Absorptivity \(\eta\)
        & \(1.0\) \\

        & Gaussian source width \(\sigma\)
        & \(50~\mu\mathrm{m}\) \\

        & Scan speed \(V\)
        & \(2~\mathrm{m\,s^{-1}}\) \\
        \hline
    \end{tabular}
\end{table}

All simulations were performed on an Intel Core i9-14900KF
CPU with 32~GB of RAM. The difference between discrete temperature fields is reported using the
relative 2-norm difference,
\begin{equation}
    e_{\mathrm{rel}}(A,B)
    =
    \frac{
        \left\|
        \mathbf{T}_{A}
        -
        \mathbf{T}_{B}
        \right\|_2
    }{
        \left\|
        \mathbf{T}_{B}
        \right\|_2
    }
    \times 100\%,
    \label{eq:relative_2norm_difference}
\end{equation}
where \(\mathbf{T}_{A}\) and \(\mathbf{T}_{B}\) are the discrete
temperature vectors associated with the two solutions being compared,
and \(B\) is the reference solution. The choice of solution \(A\) and reference solution \(B\) is specified for each comparison below.

\subsection{Numerical validation against full finite-element analysis}
\label{sec:abaqus_validation}

\begin{figure}[htbp]
    \centering

    \begin{subfigure}[t]{0.49\textwidth}
        \centering
        \includegraphics[
            width=\linewidth,
            height=0.29\textheight,
            keepaspectratio
        ]{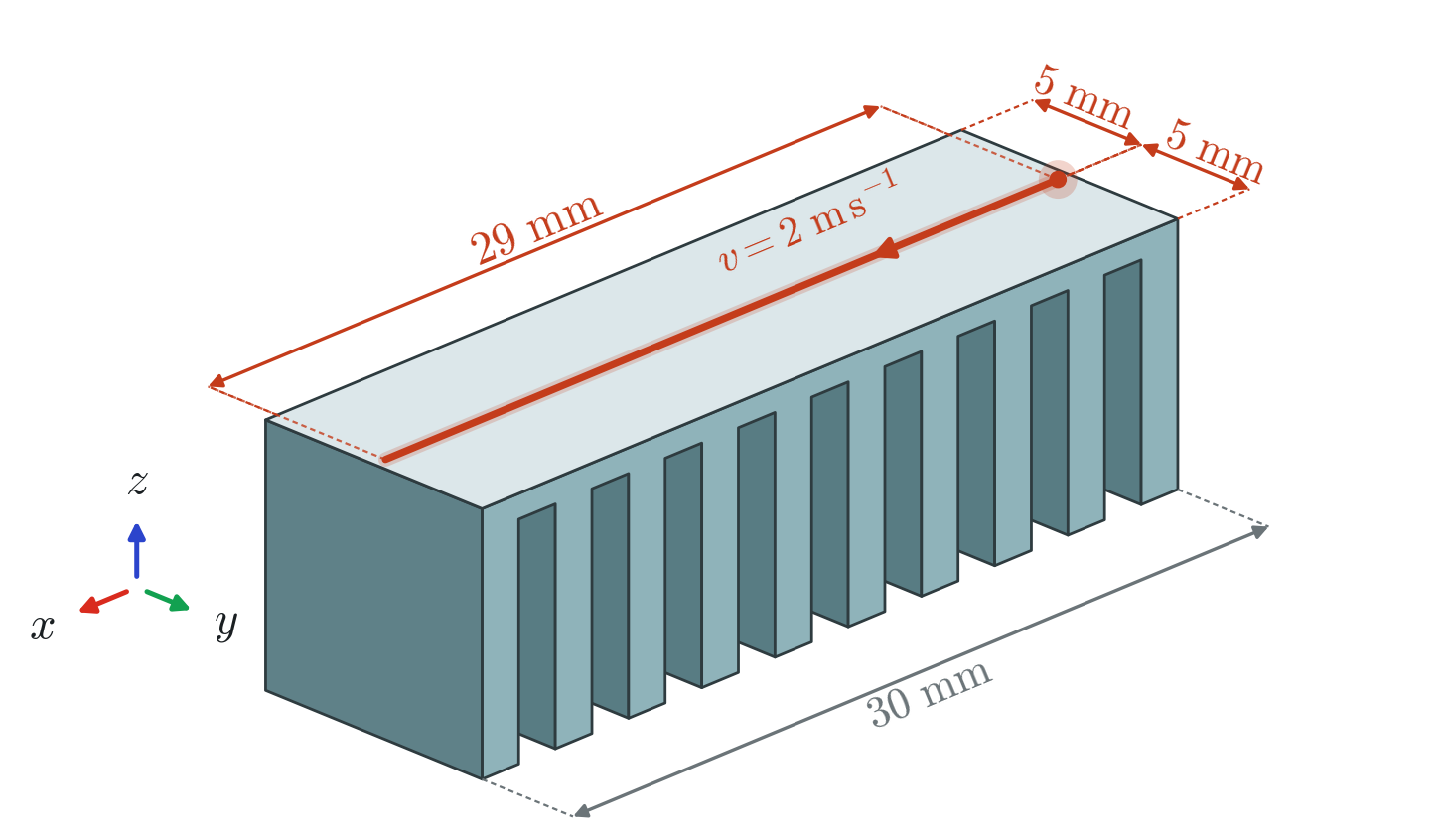}
        \caption{}
        \label{fig:abaqus_validation_setup_a}
    \end{subfigure}
    \hfill
    \begin{subfigure}[t]{0.49\textwidth}
        \centering
        \includegraphics[
            width=\linewidth,
            height=0.29\textheight,
            keepaspectratio
        ]{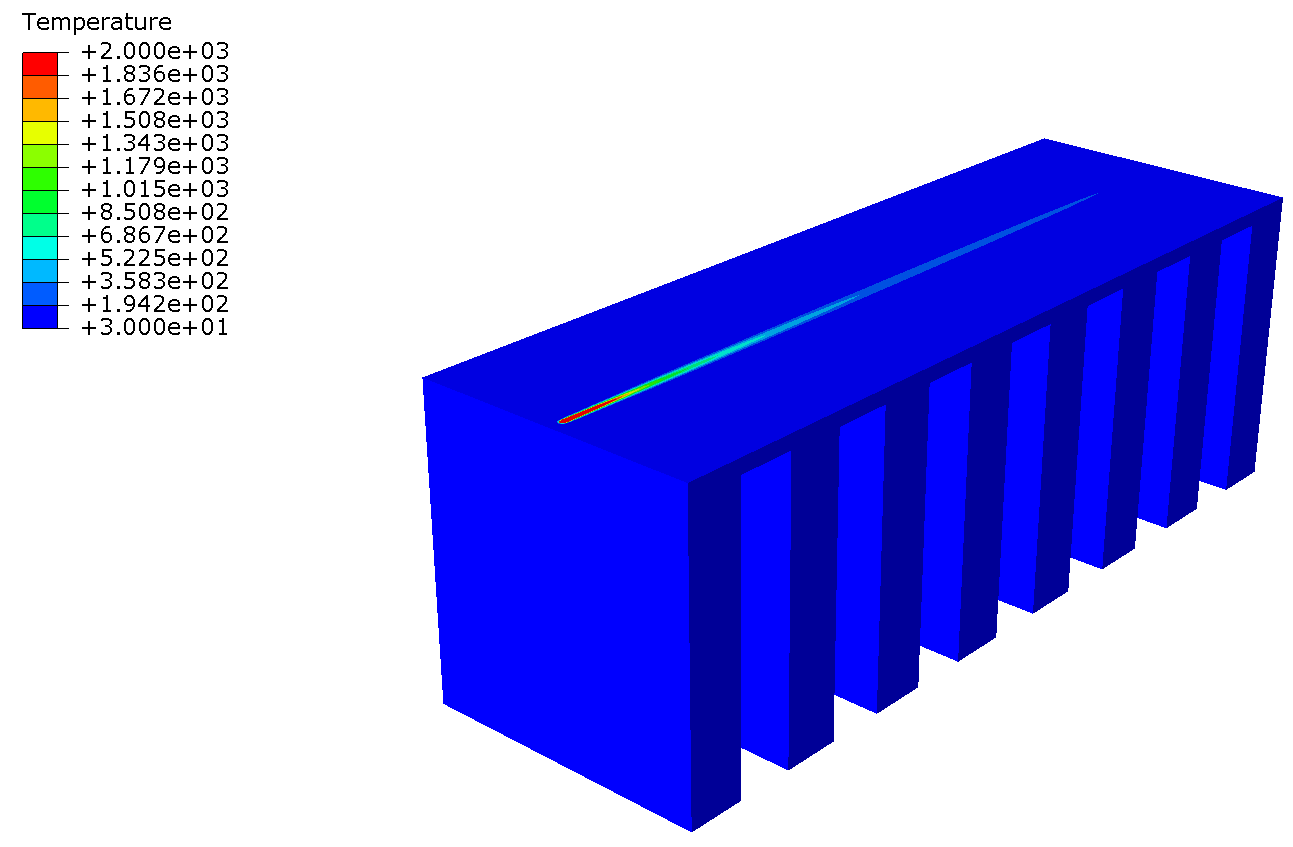}
        \caption{}
        \label{fig:abaqus_validation_setup_b}
    \end{subfigure}

    \caption{Single-line problem used for validation of the hybrid
    formulations against the full finite-element solution.
    (a) Validation geometry and scan path.
    (b) Representative Abaqus temperature field during scanning,
    displayed over \(30\)--\(2000~^{\circ}\mathrm{C}\), with higher
    temperatures shown in the upper contour color. Panel (a) was generated with assistance from Claude Opus 5 (Anthropic) and reviewed by the authors.}
    \label{fig:abaqus_validation_setup}
\end{figure}

To independently assess the temperature predictions of the hybrid
formulations, a single-scan-line problem is compared with a full
finite-element solution obtained using Abaqus.  The setup is illustrated in Fig.~\ref{fig:abaqus_validation_setup}. A
single scan line is applied to the upper surface of the fin. The laser
travels along a \(29~\mathrm{mm}\) line at \(2~\mathrm{m\,s^{-1}}\),
corresponding to a scanning duration of \(0.0145~\mathrm{s}\), followed
by \(0.2~\mathrm{s}\) of cooling. The hybrid models represent the laser through the Eagar--Tsai solution, while Abaqus applies the corresponding moving Gaussian surface flux directly.

The Abaqus model uses linear tetrahedral heat-transfer elements with
a prescribed size of \(50~\mu\mathrm{m}\) near the upper surface,
where the Gaussian heat flux is applied, and \(250~\mu\mathrm{m}\)
away from the surface. The resulting model contains
\(1{,}373{,}437\) nodes and \(7{,}171{,}717\) elements. Fixed time
increments of \(1.25\times10^{-5}~\mathrm{s}\) during scanning and
\(5\times10^{-4}~\mathrm{s}\) during cooling are used. The hybrid
simulations use a coarser linear tetrahedral mesh with prescribed
element sizes of \(200~\mu\mathrm{m}\) near the upper surface and
\(1000~\mu\mathrm{m}\) away from it. This mesh contains \(49{,}641\)
nodes and \(238{,}833\) elements, and the hybrid time-step size is
\(1.25\times10^{-4}~\mathrm{s}\).

Spatial and temporal refinement checks were performed for the Abaqus
model. Refining the prescribed element sizes from
\(100/500~\mu\mathrm{m}\) to \(50/250~\mu\mathrm{m}\) changed the
temperature field by \(0.605\%\) in the relative 2-norm, with a
root-mean-square error (RMSE) of \(0.198~^{\circ}\mathrm{C}\).
On the \(100/500~\mu\mathrm{m}\) mesh, halving the scanning and
cooling time increments to the values stated above changed the
field by \(0.0463\%\), with an RMSE of
\(0.015~^{\circ}\mathrm{C}\). These results show small differences
under the tested refinements, with greater sensitivity to spatial
than temporal discretization. The \(50/250~\mu\mathrm{m}\) solution
is therefore used as the finest tested Abaqus numerical reference.

Since the Abaqus and hybrid models use different finite-element
meshes, their temperature fields are compared at \(14{,}000\) common
spatial sampling points and 70 common saved times. The Abaqus
temperature is obtained by finite-element interpolation at the sampling
points. For the hybrid solutions, the analytical temperature is
evaluated directly at the same points and combined with the interpolated
finite-element correction. In addition to the relative 2-norm
difference defined in Eq.~\eqref{eq:relative_2norm_difference}, the
root-mean-square error (RMSE) is reported over the same temperature
samples.

\begin{table}[htbp]
    \centering
    \caption{Temperature-field comparison for the single-line
    validation problem.}
    \label{tab:abaqus_validation}
    \begin{tabular}{lcc}
        \hline
        Comparison
        & Relative 2-norm difference (\%)
        & RMSE (\(^{\circ}\mathrm{C}\)) \\
        \hline
        Conventional hybrid vs.\ Abaqus
        & \(0.330\)
        & \(0.108\) \\

        Localized hybrid vs.\ Abaqus
        & \(0.332\)
        & \(0.108\) \\

        Localized vs.\ conventional hybrid
        & \(0.0258\)
        & \(0.008\) \\
        \hline
    \end{tabular}
\end{table}

As shown in Table~\ref{tab:abaqus_validation}, the conventional and
localized hybrid formulations show close agreement with the
full finite-element solution. Their relative 2-norm differences are
\(0.330\%\) and \(0.332\%\), respectively, with RMSE values of
approximately \(0.108~^{\circ}\mathrm{C}\).

The difference between the localized and conventional hybrid solutions
is considerably smaller. The relative 2-norm difference is
\(0.0258\%\), with an RMSE of \(0.008~^{\circ}\mathrm{C}\). These
results indicate that localization introduces only a small additional
difference relative to the conventional hybrid formulation, while both
hybrid formulations remain in close agreement with the independent
full finite-element solution.

The Abaqus simulation required approximately \(88{,}848~\mathrm{s}\)
(\(24.7~\mathrm{h}\)) of wall-clock time, compared with
\(1576.2~\mathrm{s}\) for the conventional hybrid simulation and
\(1281.6~\mathrm{s}\) for the localized hybrid simulation. The full finite-element simulation required
substantially more computational time. This reflects the cost of
directly resolving the localized moving heat source with the finer
spatial and temporal discretization required by the full finite-element
model.

\subsection{Comparison with temporal truncation}
\label{sec:temporal_truncation}

\begin{figure}[htbp]
    \centering

    \begin{subfigure}[c]{0.48\linewidth}
        \centering
        \includegraphics[
            width=\linewidth,
            height=0.30\textheight,
            keepaspectratio
        ]{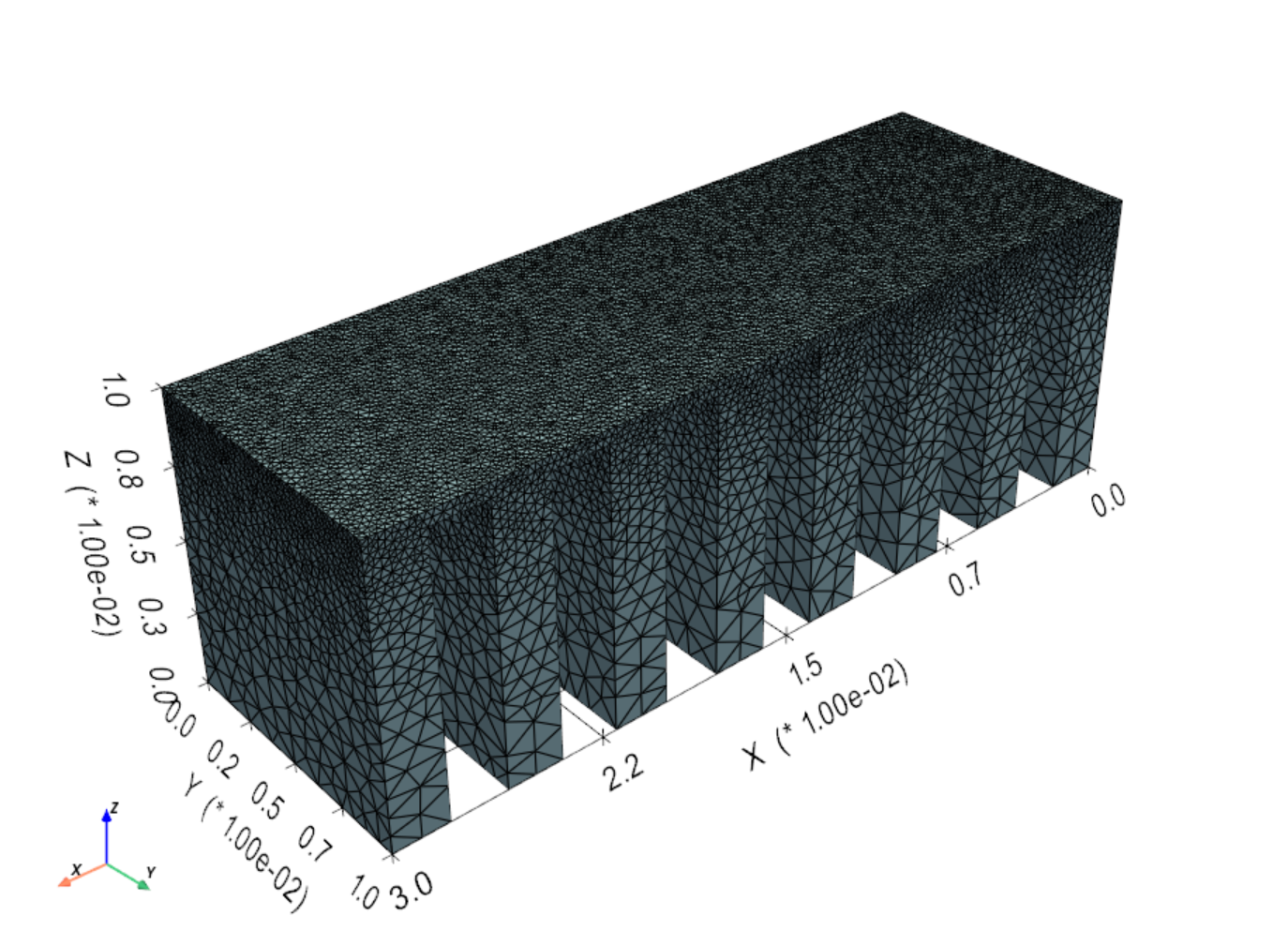}
        \caption{}
        \label{fig:problem_setup_mesh}
    \end{subfigure}
    \hfill
    \begin{subfigure}[c]{0.48\linewidth}
        \centering
        \includegraphics[
            width=\linewidth,
            height=0.30\textheight,
            keepaspectratio
        ]{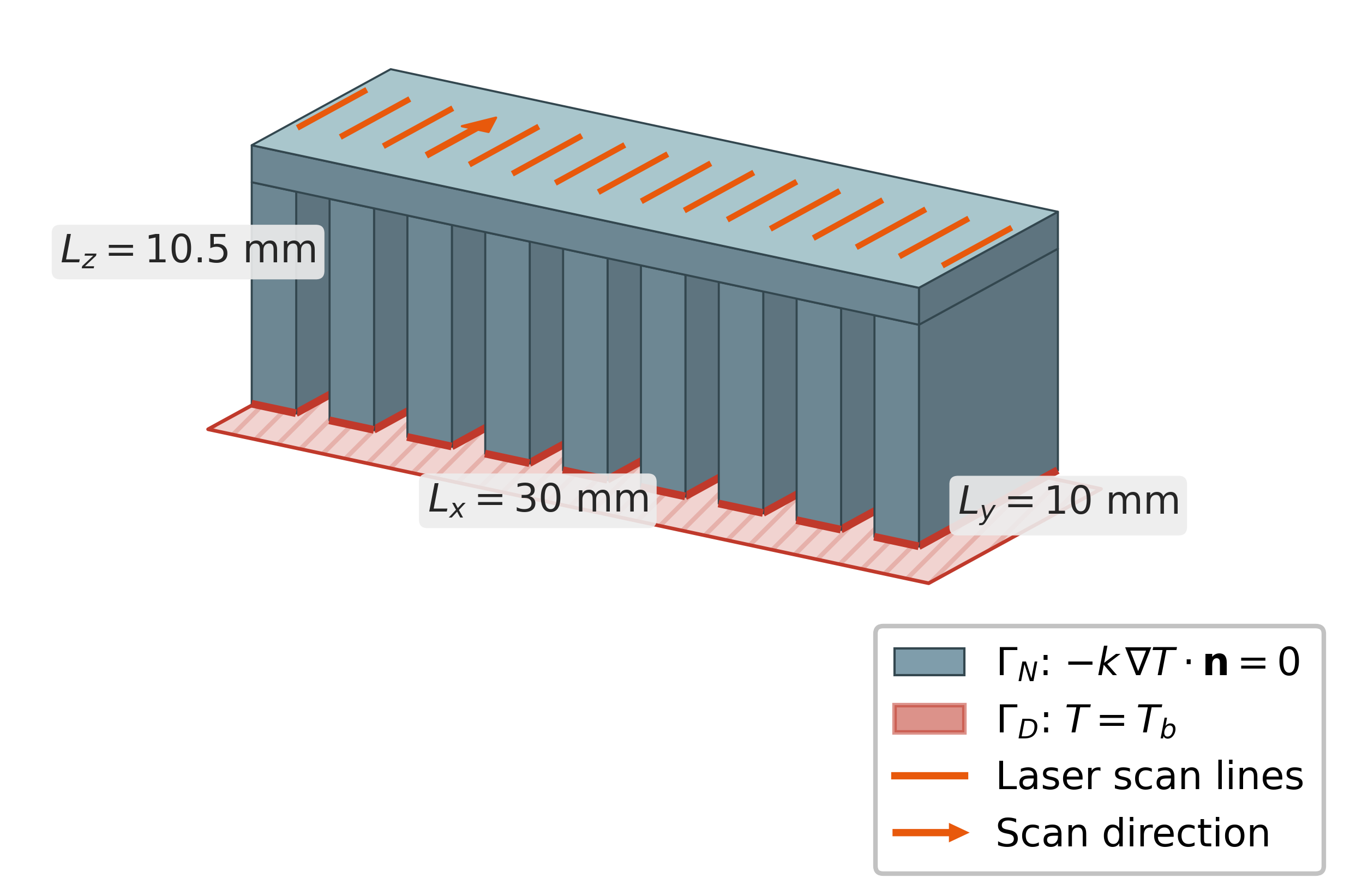}
        \caption{}
        \label{fig:problem_setup_geometry}
    \end{subfigure}

    \caption{Fin geometry, finite-element mesh, and scan arrangement for the multiline studies. (a) Linear tetrahedral mesh with finer resolution near the upper surface and coarser resolution toward the bottom. (b) Schematic scan-line arrangement, scan direction, and
    thermal boundary conditions. Panel (b) was generated with assistance from Claude Opus 5 (Anthropic) and reviewed by the authors.}
    \label{fig:problem_setup}
\end{figure}

This study and the studies in
Sections~\ref{sec:localization_parameter_study} and
\ref{sec:detailed_500_line_case} use the parallel scan-line
arrangement and finite-element mesh shown in
Fig.~\ref{fig:problem_setup}, with the numerical settings listed
in Table~\ref{tab:computational_parameters}.
All three studies use 500 scan lines and include
\(1~\mathrm{s}\) of post-scan cooling. Each line extends in the
\(y\)-direction from \(y=1~\mathrm{mm}\) to \(y=9~\mathrm{mm}\),
leaving \(1~\mathrm{mm}\) margins at both ends. The lines are
uniformly spaced in the \(x\)-direction between
\(x=0.5~\mathrm{mm}\) and \(x=29.5~\mathrm{mm}\), leaving
\(0.5~\mathrm{mm}\) margins at the corresponding edges.
Together, they span a \(29~\mathrm{mm}\times8~\mathrm{mm}\)
region of the \(30~\mathrm{mm}\times10~\mathrm{mm}\) upper surface.

\begin{table}[htbp]
    \centering
    \caption{Numerical settings shared by the multiline hybrid studies.}
    \label{tab:computational_parameters}
    \begin{tabular}{lll}
        \hline
        Category & Parameter & Value \\
        \hline
        Scan path
        & Line length
        & \(8~\mathrm{mm}\) \\

        & Hatch spacing
        & \(58~\mu\mathrm{m}\) \\

        & Delay between lines
        & \(0.5~\mathrm{ms}\) \\
        \hline
        Finite-element model
        & Prescribed element-size range
        & \(200\text{--}1000~\mu\mathrm{m}\) \\

        & Number of nodes
        & \(42470\) \\

        & Number of tetrahedral elements
        & \(202094\) \\
        \hline
        Time integration
        & Time-step size
        & \(0.5~\mathrm{ms}\) \\
        \hline
    \end{tabular}
\end{table}

For comparisons between the localized and conventional hybrid
solutions, the relative 2-norm difference in
Eq.~\eqref{eq:relative_2norm_difference} is evaluated using nodal
temperatures over all saved times, with the conventional solution as
the reference.

Computational performance is quantified using the analytical speedup
\begin{equation}
    S_{\mathrm{A}}
    =
    \frac{
        t_{\mathrm{A,conv}}
    }{
        t_{\mathrm{A,loc}}
    },
    \label{eq:analytical_speedup}
\end{equation}
and the overall speedup
\begin{equation}
    S_{\mathrm{overall}}
    =
    \frac{
        t_{\mathrm{conv}}
    }{
        t_{\mathrm{loc}}
    },
    \label{eq:overall_speedup}
\end{equation}
where \(t_{\mathrm{A,conv}}\) and \(t_{\mathrm{A,loc}}\) are the times
spent evaluating analytical temperatures and gradients, and
\(t_{\mathrm{conv}}\) and \(t_{\mathrm{loc}}\) are the total
computational times of the conventional and localized simulations,
respectively. Analytical work is also
measured by the number of line--point evaluations, where each
evaluation corresponds to the contribution of one scan line at one
spatial point. The count is accumulated over all time steps and
includes mesh nodes and the boundary and volume integration points
required by the hybrid formulation.

The first multiline study compares the localized formulation
with uncompensated temporal truncation for the 500-line problem.
The comparison examines whether discarding older analytical
contributions without accounting for their removal can maintain
agreement with the conventional hybrid solution. The localized
formulation uses a spatial core radius and transition width of
\(1~\mathrm{mm}\) each, and a temporal retention age and transfer
duration of \(50~\mathrm{ms}\) each. The sensitivity to these
parameters and the reason for choosing these values are examined
in Section~\ref{sec:localization_parameter_study}.

In the temporal-truncation formulation, each line is evaluated throughout the computational domain until a cutoff interval
\(\tau_c=t_r+\Delta t_\chi=0.1~\mathrm{s}\) after its completion.
This interval matches the selected localized formulation, which
uses \(50~\mathrm{ms}\) of retention followed by
\(50~\mathrm{ms}\) of gradual transfer. The retained analytical
field is
\begin{equation}
    T_A^{\mathrm{trunc}}(\mathbf{x},t)
    =
    \sum_{\substack{
        \ell\in\mathcal{A}(t)\\
        t<t_{e,\ell}+\tau_c
    }}
    T_\ell(\mathbf{x},t).
    \label{eq:truncated_analytical_field}
\end{equation}
The finite-element correction is advanced using
Eq.~\eqref{eq:conventional_FE_equation}, with the boundary
conditions in Eqs.~\eqref{eq:conventional_dirichlet} and
\eqref{eq:conventional_neumann} evaluated using
\(T_A^{\mathrm{trunc}}\). No localization residual or compensating
transfer is applied when a line is removed. No spatial localization
is used. Thus, both methods use the same removal deadline, but localization
accounts for the gradual reduction of each analytical contribution.

A completed scan line continues to affect the temperature field
as its heat diffuses through the domain. Deleting its analytical
contribution without transferring the removed field to the
finite-element correction therefore loses part of the thermal
history. Updating the boundary correction does not replace the
missing contribution throughout the interior of the domain.
Repeated deletions can therefore accumulate into a substantial
discrepancy from the conventional solution.
In the localized formulation, the temporal residual in
Eq.~\eqref{eq:temporal_residual} accounts for the decreasing
analytical weight, together with the spatial residual and
boundary conditions derived in Section~\ref{sec:methodology}.

\begin{figure}[htbp]
    \centering
    \includegraphics[width=0.85\linewidth]
    {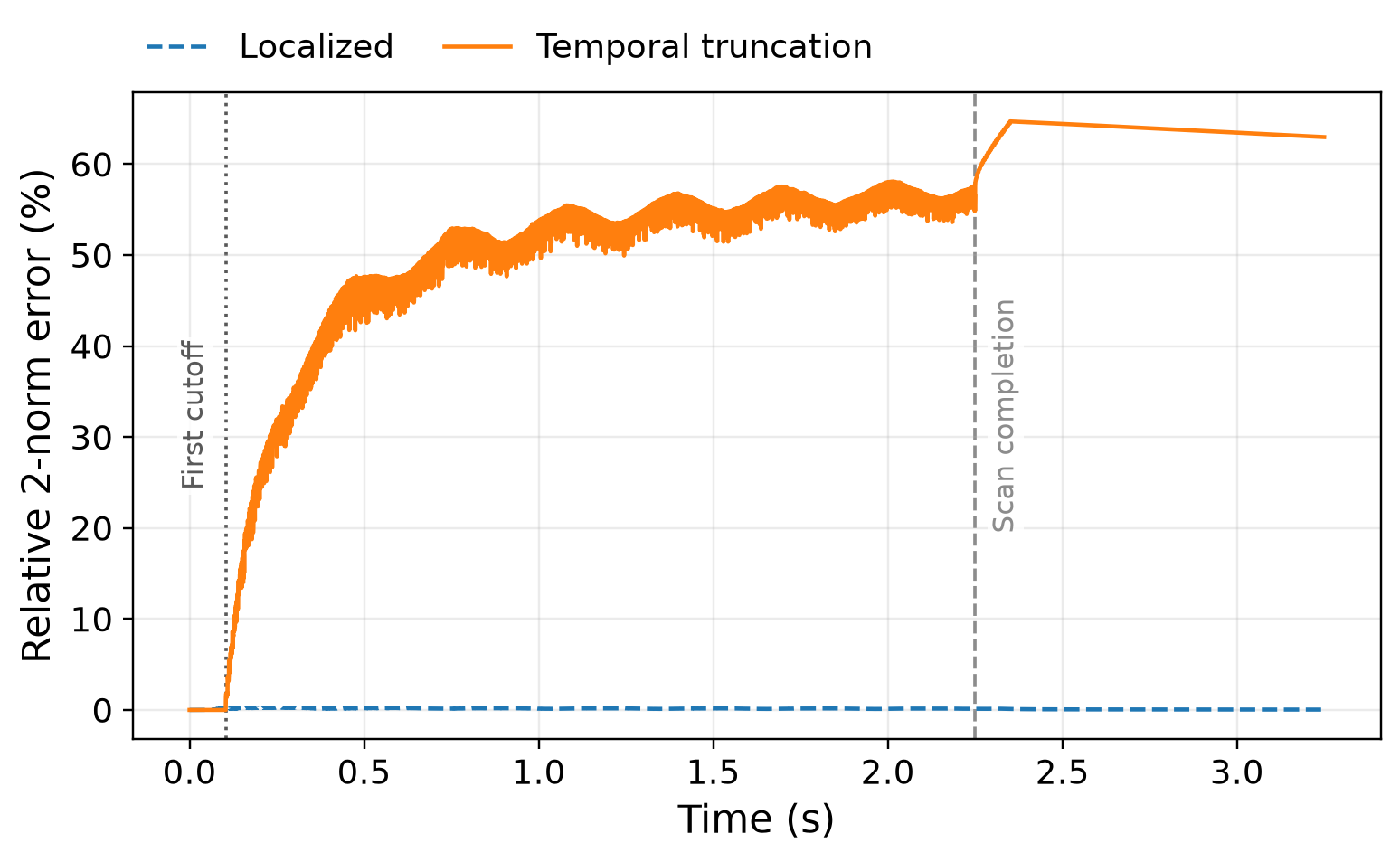}
    \caption{Relative nodal 2-norm difference from the conventional
    hybrid solution at each saved time for the 500-line case.
    The vertical lines indicate the first truncation event and
    scan completion.}
    \label{fig:temporal_truncation_error}
\end{figure}

Figure~\ref{fig:temporal_truncation_error} shows the relative
2-norm difference at each saved time. Before the first deletion at
\(t=0.104~\mathrm{s}\), temporal truncation agrees with the
conventional solution to numerical precision. The difference then
increases as scan lines are removed, reaching approximately
\(47.2\%\) at \(t=0.5~\mathrm{s}\) and \(54.9\%\) at scan
completion, \(t=2.2495~\mathrm{s}\). The localized solution
remains close to the conventional solution throughout the simulation.

This demonstrates the importance of compensating for
removed analytical contributions when limiting the retained
scan history.

\subsection{Effect of localization parameters}
\label{sec:localization_parameter_study}

This study examines how the spatial and temporal localization
parameters in Algorithm~\ref{alg:localized_hybrid} affect accuracy
and computational cost. The localization parameters are varied
while the scan path, material properties, mesh, and numerical
integration settings are held fixed.

A common scale factor \(\lambda\) is applied to the four localization
parameters. At \(\lambda=1\), the spatial core radius and transition
width are both \(1~\mathrm{mm}\), while the temporal retention age
and transfer duration are both \(50~\mathrm{ms}\). Each parameter
is multiplied by \(\lambda\), with
\[
    \lambda\in\{0.5,\;0.75,\;1.0,\;1.25,\;1.5,\;1.75,\;2.0\}.
\]
Increasing \(\lambda\) enlarges the spatial region and extends
the time interval over which analytical contributions are evaluated.
Larger spatial regions allow each contribution to diffuse further
before it is represented through the finite-element correction.
Longer temporal transition intervals also reduce the rate at which
its analytical weight changes, as reflected in the temporal residual
in Eq.~\eqref{eq:temporal_residual}. For a fixed mesh and time step,
these changes are expected to improve agreement with the conventional
solution, while requiring more analytical evaluation. All localized
results are compared with the same conventional hybrid reference.

\begin{figure}[htbp]
    \centering
    \begin{subfigure}[t]{0.49\linewidth}
        \centering
        \includegraphics[width=\linewidth]
        {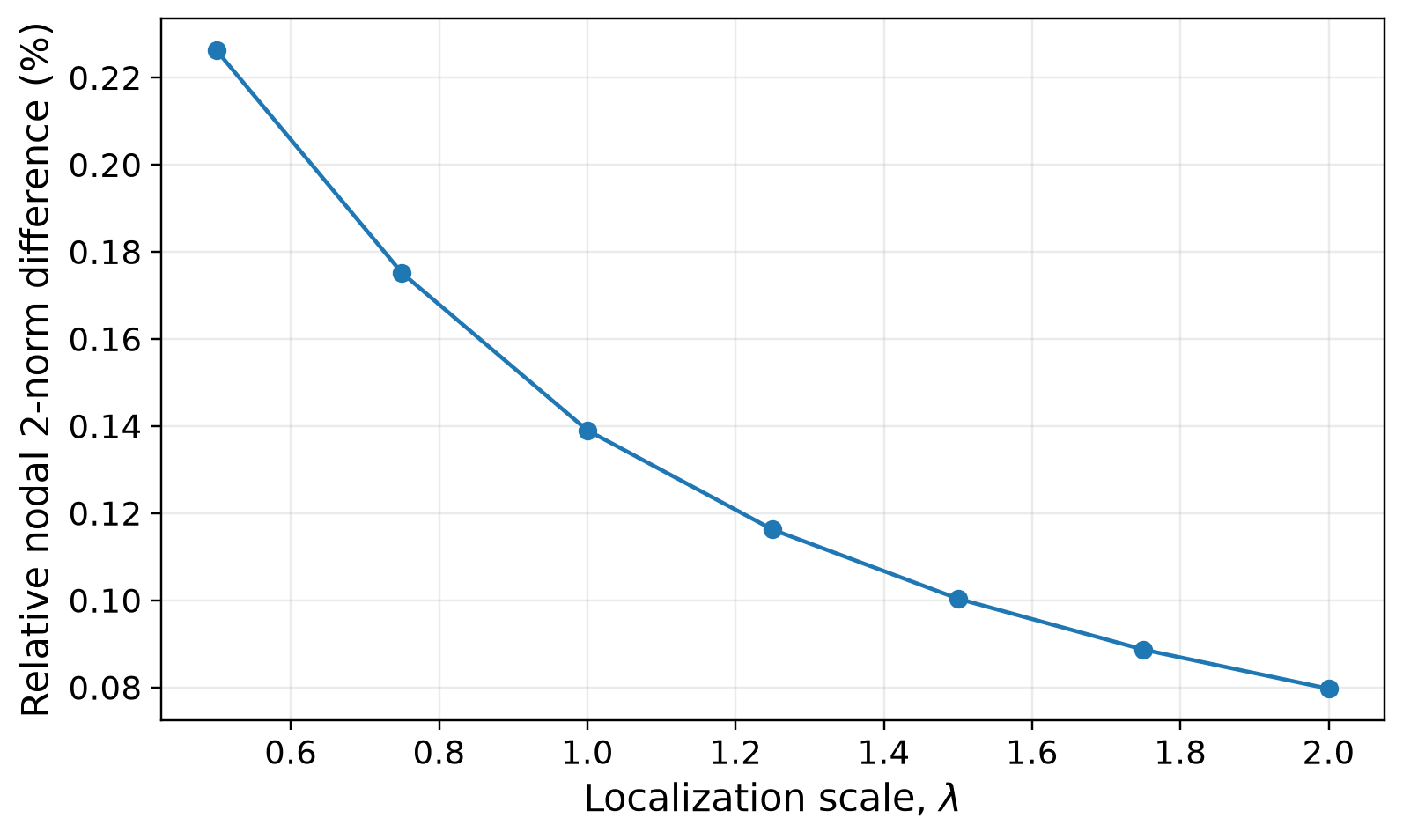}
        \caption{}
        \label{fig:localization_parameter_difference}
    \end{subfigure}
    \hfill
    \begin{subfigure}[t]{0.49\linewidth}
        \centering
        \includegraphics[width=\linewidth]
        {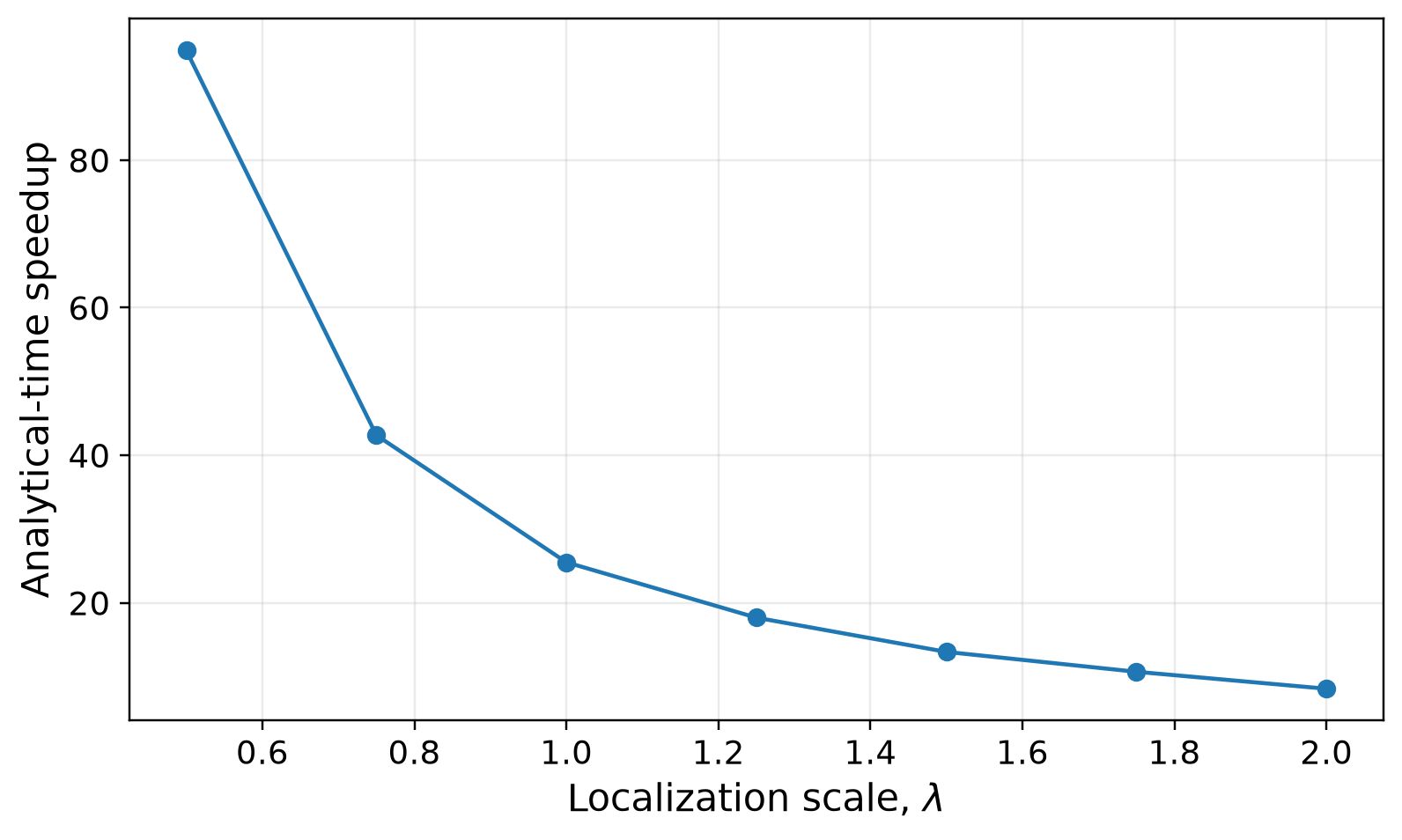}
        \caption{}
        \label{fig:localization_parameter_speedup}
    \end{subfigure}

    \caption{Influence of the localization scale \(\lambda\) in the
    500-line case: (a) relative nodal 2-norm difference from the
    conventional hybrid solution over all saved times and
    (b) analytical speedup relative to the conventional formulation.}
    \label{fig:localization_parameter_study}
\end{figure}

Figure~\ref{fig:localization_parameter_study} shows how accuracy
and analytical speedup change with the localization scale.
The relative nodal 2-norm difference, evaluated over all saved
times, decreases monotonically from \(0.226\%\) at
\(\lambda=0.5\) to \(0.0797\%\) at \(\lambda=2\).
Over the same range, the number of analytical line--point
evaluations increases from \(4.65\times10^{8}\) to
\(1.03\times10^{10}\), and the analytical speedup decreases
from \(94.8\) to \(8.37\). These results are consistent with
the expected tradeoff: larger localization scales give closer
agreement at the cost of more analytical evaluation.
Substantial analytical savings remain throughout the tested range.

The setting \(\lambda=1\) is used in the temporal-truncation
comparison in Section~\ref{sec:temporal_truncation} and the
computational-savings study in
Section~\ref{sec:detailed_500_line_case}.
Although \(\lambda=0.5\) provides a much larger analytical
speedup, its overall speedup (when also accounting for FEA) increases only from \(3.92\)
to \(4.59\). Thus, \(\lambda=1\) provides closer agreement
at a modest additional computational cost.

\subsection{Analytical history and computational savings}
\label{sec:detailed_500_line_case}

The same 500-line case used in Section~\ref{sec:temporal_truncation} is examined here to show how localization limits the analytical history and
reduces computational cost. The localization parameters are fixed at the \(\lambda=1\) values discussed in Section~\ref{sec:localization_parameter_study}.

\begin{figure}[htbp]
    \centering
    \includegraphics[width=0.72\linewidth]
    {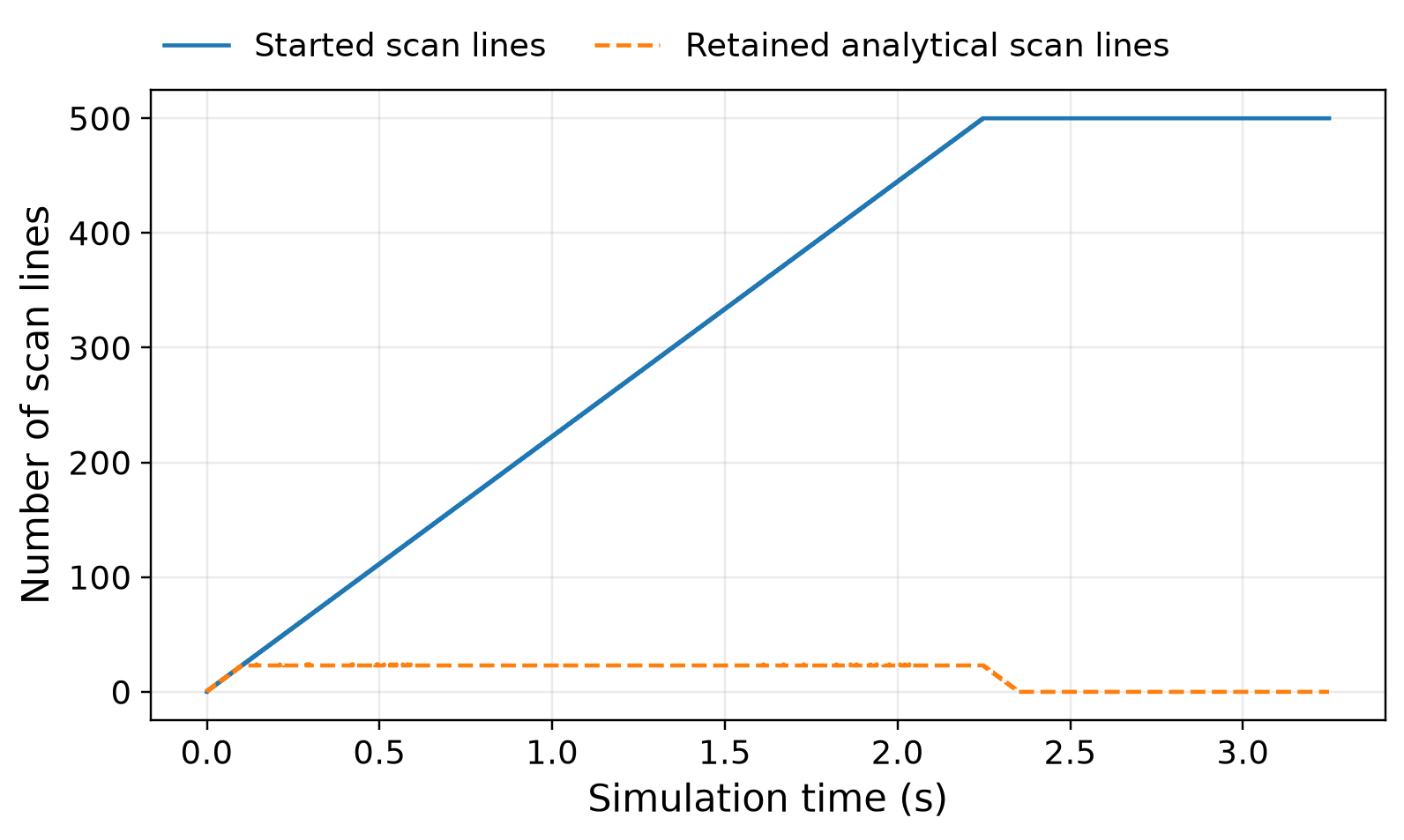}
    \caption{Number of started scan lines and explicitly retained
    analytical scan lines for the localized formulation in the
    500-line case.}
    \label{fig:retained_line_history}
\end{figure}

Figure~\ref{fig:retained_line_history} compares the number of scan
lines that have begun with the number that remain explicitly included
in the localized analytical field. In the conventional hybrid
formulation, every started line remains in the analytical superposition
throughout the simulation. The number of analytical contributions
therefore increases to 500 by the end of scanning. In the localized
formulation, completed lines are gradually removed from the analytical
field, and the number of retained analytical lines never exceeds 24.

Table~\ref{tab:detailed_timing_500} compares the time spent in the
main computational components of the two formulations. Localization
reduces the analytical evaluation time from \(13059.6~\mathrm{s}\)
to \(512.1~\mathrm{s}\), giving an analytical speedup of \(25.5\).
The finite-element solution times remain similar, at
\(3308.9~\mathrm{s}\) and \(3209.0~\mathrm{s}\), respectively.

\begin{table}[htbp]
    \centering
    \caption{Time spent in the main computational components for the
    500-line case. Localization and assembly includes localization
    operations, localization-residual assembly, and boundary assembly.}
    \label{tab:detailed_timing_500}
    \begin{tabular}{lrr}
        \hline
        Component
        & Conventional (s)
        & Localized (s) \\
        \hline
        Analytical evaluation
        & 13059.6 & 512.1 \\
        Finite-element solution
        & 3308.9 & 3209.0 \\
        Localization and assembly
        & 4.8 & 280.0 \\
        \hline
    \end{tabular}
\end{table}

Localization introduces additional work to construct the localization
functions and assemble the residual terms. The combined localization
and assembly time is \(280.0~\mathrm{s}\), compared with
\(4.8~\mathrm{s}\) for boundary assembly in the conventional
formulation. The total simulation time decreases from
\(16375.5~\mathrm{s}\) to \(4176.4~\mathrm{s}\), giving an overall
speedup of \(3.92\). Over all saved times
and mesh nodes, the relative nodal 2-norm difference is
\(0.139\%\).

\section{Conclusions}
\label{sec:conclusions}

This work developed a space--time localized hybrid
analytical--finite-element formulation for part-scale thermal simulation of
laser powder bed fusion. The formulation addresses the growing
analytical workload of conventional hybrid methods as scan-line
contributions accumulate. Spatial and temporal localization
restrict where and how long each contribution is evaluated
analytically. Residual terms and boundary conditions derived
from the heat equation allow the finite-element correction to
account for the thermal effect of contributions removed from
the analytical field.

The proposed formulation maintained close agreement with an
independent full finite-element solution for the single-line
benchmark. Both the conventional and localized hybrid
formulations reproduced the reference temperature field closely,
supporting the use of the localized formulation in the subsequent
multiline comparisons.

Limiting analytical history requires accounting for the thermal
effect of removed contributions. Uncompensated temporal truncation
produced a large discrepancy from the conventional hybrid solution,
whereas localization maintained close agreement. These results
demonstrate the importance of the residual correction in reducing
the explicit analytical history while retaining its thermal effect.

The extent of localization controls the tradeoff between agreement
and computational cost. Larger spatial regions and longer temporal
intervals improved agreement with the conventional solution,
while smaller settings provided greater analytical savings.
Substantial savings were obtained throughout the tested parameter
range, allowing the parameters to be chosen according to the
desired balance between agreement and computational effort.

With the selected parameters, the retained analytical history
never exceeded 24 line contributions in the 500-line case,
while the conventional formulation retained all 500 contributions
by the end of scanning. This reduction
gave an analytical speedup of \(25.5\) and an overall speedup
of \(3.92\), with a relative nodal 2-norm difference of
\(0.139\%\). The additional localization and assembly work
was small compared with the analytical evaluation time saved.
The finite-element solution time remained nearly unchanged
and became the dominant computational cost, limiting the
overall speedup. This provides a basis for applying
the proposed method to larger part models (with a much larger number of scan lines), with greater analytical
savings expected as the conventional history grows.

Possible future directions include multilayer builds with evolving
geometry, more varied scan paths, and adaptive selection of the
localization parameters. Adjusting the spatial regions and transfer
intervals to the thermal field, mesh, and time step could provide
further savings while controlling the difference from the
conventional solution. Extending the formulation to
temperature-dependent material properties and phase change would
broaden its physical scope, while experimental validation would
help assess its predictions under practical LPBF conditions.

\subsection*{Acknowledgments}
The authors acknowledge support from the US Office of Naval Research under Grant N000142412391.

\subsection*{Conflicts of Interest}
The authors declare no conflicts of interest.

\subsection*{Data Availability Statement}
The code developed in this work can be made available upon request.

\subsection*{Use of AI Tools}

The authors used generative AI assistants, specifically, ChatGPT 6 (OpenAI) and Claude Opus 5 (Anthropic), during the preparation of this manuscript. The tools were used for code development, drafting of text, and generating schematic illustrations. All content was reviewed, verified, and edited by the authors, who take full responsibility for the accuracy and integrity of the work.

\bibliographystyle{elsarticle-num}
\bibliography{references}

@article{subedi2022towards,
  title={Towards the optimal design of support structures for laser powder bed fusion-based metal additive manufacturing via thermal equivalent static loads},
  author={Subedi, Subodh C and Shahba, Ahmad and Thevamaran, Mythili and Thoma, Dan J and Suresh, Krishnan},
  journal={Additive Manufacturing},
  volume={57},
  pages={102956},
  year={2022},
  publisher={Elsevier}
}

@article{prabhune2020fast,
  title={A fast matrix-free elasto-plastic solver for predicting residual stresses in additive manufacturing},
  author={Prabhune, Bhagyashree C and Suresh, Krishnan},
  journal={Computer-Aided Design},
  volume={123},
  pages={102829},
  year={2020},
  publisher={Elsevier}
}

@article{subedi2024deformation,
  author  = {Subedi, Subodh C. and Thoma, Dan J. and Suresh, Krishnan},
  title   = {Deformation constrained support-structure optimization for laser powder bed fusion},
  journal = {Additive Manufacturing},
  volume  = {89},
  pages   = {104294},
  year    = {2024},
  doi     = {10.1016/j.addma.2024.104294}
}

@article{busachi2017review,
  title={A review of Additive Manufacturing technology and Cost Estimation techniques for the defence sector},
  author={Busachi, Alessandro and Erkoyuncu, John and Colegrove, Paul and Martina, Filomeno and Watts, Chris and Drake, Richard},
  journal={CIRP Journal of Manufacturing Science and Technology},
  volume={19},
  pages={117--128},
  year={2017},
  publisher={Elsevier}
}

@incollection{kumar2016current,
  title={Current trends of additive manufacturing in the aerospace industry},
  author={Kumar, L Jyothish and Krishnadas Nair, CG},
  booktitle={Advances in 3D printing \& additive manufacturing technologies},
  pages={39--54},
  year={2016},
  publisher={Springer}
}

@article{pereira2019comparison,
  title={A comparison of traditional manufacturing vs additive manufacturing, the best method for the job},
  author={Pereira, Tanisha and Kennedy, John V and Potgieter, Johan},
  journal={Procedia manufacturing},
  volume={30},
  pages={11--18},
  year={2019},
  publisher={Elsevier}
}

@article{promoppatum2017comprehensive,
  title={A comprehensive comparison of the analytical and numerical prediction of the thermal history and solidification microstructure of Inconel 718 products made by laser powder-bed fusion},
  author={Promoppatum, Patcharapit and Yao, Shi-Chune and Pistorius, P Chris and Rollett, Anthony D},
  journal={Engineering},
  volume={3},
  number={5},
  pages={685--694},
  year={2017},
  publisher={Elsevier}
}

@article{soundararajan2021review,
  title={Review on modeling techniques for powder bed fusion processes based on physical principles},
  author={Soundararajan, Balaji and Sofia, Daniele and Barletta, Diego and Poletto, Massimo},
  journal={Additive Manufacturing},
  volume={47},
  pages={102336},
  year={2021},
  publisher={Elsevier}
}

@article{gong2014analysis,
  title={Analysis of defect generation in Ti--6Al--4V parts made using powder bed fusion additive manufacturing processes},
  author={Gong, Haijun and Rafi, Khalid and Gu, Hengfeng and Starr, Thomas and Stucker, Brent},
  journal={Additive Manufacturing},
  volume={1},
  pages={87--98},
  year={2014},
  publisher={Elsevier}
}

@article{scheel2020adaptive,
  title={Adaptive local-global multiscale approach for thermal simulation of the selective laser melting process},
  author={Scheel, P and Mazza, Edoardo and Hosseini, Ehsan},
  journal={Additive Manufacturing},
  volume={36},
  pages={101518},
  year={2020},
  publisher={Elsevier}
}

@article{bresson2023global,
  title={Global-to-local simulation of the thermal history in the laser powder bed fusion process based on a multiscale finite element approach},
  author={Bresson, Yves and Tongne, Am{\`e}vi and Baili, Maher and Arnaud, Lionel},
  journal={The International Journal of Advanced Manufacturing Technology},
  volume={127},
  number={9},
  pages={4727--4744},
  year={2023},
  publisher={Springer}
}

@article{li2016multiscale,
  title={A multiscale modeling approach for fast prediction of part distortion in selective laser melting},
  author={Li, Chao and Fu, CH and Guo, YB and Fang, FZ},
  journal={Journal of materials processing technology},
  volume={229},
  pages={703--712},
  year={2016},
  publisher={Elsevier}
}

@article{xia2017porosity,
  title={Porosity evolution and its thermodynamic mechanism of randomly packed powder-bed during selective laser melting of Inconel 718 alloy},
  author={Xia, Mujian and Gu, Dongdong and Yu, Guanqun and Dai, Donghua and Chen, Hongyu and Shi, Qimin},
  journal={International Journal of Machine Tools and Manufacture},
  volume={116},
  pages={96--106},
  year={2017},
  publisher={Elsevier}
}

@article{king2015overview,
  title={Overview of modelling and simulation of metal powder bed fusion process at Lawrence Livermore National Laboratory},
  author={King, Wayne and Anderson, Andrew T and Ferencz, Robert M and Hodge, Neil E and Kamath, Chandrika and Khairallah, Saad A},
  journal={Materials Science and Technology},
  volume={31},
  number={8},
  pages={957--968},
  year={2015},
  publisher={SAGE Publications Sage UK: London, England}
}

@article{love2025advancing,
  title={Advancing Metal Additive Manufacturing: A Review of Numerical Methods in DED, WAAM, and PBF},
  author={Love, Allen and Valdez Pastrana, Omar Alejandro and Behseresht, Saeed and Park, Young Ho},
  journal={Metrology},
  volume={5},
  number={2},
  pages={30},
  year={2025},
  publisher={MDPI}
}

@article{bayat2021review,
  title={A review of multi-scale and multi-physics simulations of metal additive manufacturing processes with focus on modeling strategies},
  author={Bayat, Mohamad and Dong, Wen and Thorborg, Jesper and To, Albert C and Hattel, Jesper H},
  journal={Additive Manufacturing},
  volume={47},
  pages={102278},
  year={2021},
  publisher={Elsevier}
}

@article{zagade2024analytical,
  title={Analytical modelling of scanning strategy effect on temperature field and melt track dimensions in laser powder bed fusion},
  author={Zagade, Pramod R and Gautham, BP and De, Amitava and DebRoy, Tarasankar},
  journal={Additive Manufacturing},
  volume={82},
  pages={104046},
  year={2024},
  publisher={Elsevier}
}

@article{de2025review,
  title={A Review of Modeling, Simulation, and Process Qualification of Additively Manufactured Metal Components via the Laser Powder Bed Fusion Method},
  author={De Leon, Emmanuel and Riensche, Alex and Bevans, Benjamin D and Billings, Christopher and Siddique, Zahed and Liu, Yingtao},
  journal={Journal of Manufacturing and Materials Processing},
  volume={9},
  number={1},
  pages={22},
  year={2025},
  publisher={MDPI}
}

@article{yang2018semi,
  title={A semi-analytical thermal modelling approach for selective laser melting},
  author={Yang, Y and Knol, MF and Van Keulen, F and Ayas, C},
  journal={Additive manufacturing},
  volume={21},
  pages={284--297},
  year={2018},
  publisher={Elsevier}
}

@article{moran2021scan,
  title={Scan-by-scan part-scale thermal modelling for defect prediction in metal additive manufacturing},
  author={Moran, TP and Warner, DH and Phan, N},
  journal={Additive Manufacturing},
  volume={37},
  pages={101667},
  year={2021},
  publisher={Elsevier}
}

@article{moran2018utility,
  title={Utility of superposition-based finite element approach for part-scale thermal simulation in additive manufacturing},
  author={Moran, TP and Li, P and Warner, DH and Phan, N},
  journal={Additive Manufacturing},
  volume={21},
  pages={215--219},
  year={2018},
  publisher={Elsevier}
}

@article{wang2023analytical,
  author  = {Wang, Yongfu and Ji, Xia and Liang, Steven Y.},
  title   = {Analytical Modeling of Temperature Distribution in Laser
             Powder Bed Fusion with Different Scan Strategies},
  journal = {Optics \& Laser Technology},
  volume  = {157},
  pages   = {108708},
  year    = {2023},
  doi     = {10.1016/j.optlastec.2022.108708}
}

@article{yang2020computationally,
  title={A computationally efficient thermal model for selective laser melting},
  author={Yang, Yabin and van Keulen, Fred and Ayas, Can},
  journal={Additive Manufacturing},
  volume={31},
  pages={100955},
  year={2020},
  publisher={Elsevier}
}

@article{yang2025computational,
  title={Computational efficient process simulation of geometrically complex parts in metal additive manufacturing},
  author={Yang, Y and Ji, Y and M{\"o}ller, M and Ayas, C},
  journal={International Journal of Heat and Mass Transfer},
  volume={248},
  pages={127059},
  year={2025},
  publisher={Elsevier}
}

@article{nasiri2019powerful,
  title={Powerful analytical solution to heat flow problem in welding},
  author={Nasiri, Mohammad Bagher and Enzinger, Norbert},
  journal={International Journal of Thermal Sciences},
  volume={135},
  pages={601--612},
  year={2019},
  publisher={Elsevier}
}

@article{steuben2019enriched,
  title={Enriched analytical solutions for additive manufacturing modeling and simulation},
  author={Steuben, John C and Birnbaum, Andrew J and Michopoulos, John G and Iliopoulos, Athanasios P},
  journal={Additive Manufacturing},
  volume={25},
  pages={437--447},
  year={2019},
  publisher={Elsevier}
}

@article{rosenthal1946theory,
  title={The theory of moving sources of heat and its application to metal treatments},
  author={Rosenthal, Daniel},
  journal={Transactions of the American Society of Mechanical Engineers},
  volume={68},
  number={8},
  pages={849--865},
  year={1946},
  publisher={American Society of Mechanical Engineers}
}

@article{eagar1983temperature,
  title={Temperature fields produced by traveling distributed heat sources},
  author={Eagar, TW and Tsai, NS and others},
  journal={Welding journal},
  volume={62},
  number={12},
  pages={346--355},
  year={1983}
}

@article{olleak2020part,
  title={Part-scale finite element modeling of the selective laser melting process with layer-wise adaptive remeshing for thermal history and porosity prediction},
  author={Olleak, Alaa and Xi, Zhimin},
  journal={Journal of Manufacturing Science and Engineering},
  volume={142},
  number={12},
  pages={121006},
  year={2020},
  publisher={American Society of Mechanical Engineers}
}

@article{dugast2021part,
  title={Part-scale thermal process modeling for laser powder bed fusion with matrix-free method and GPU computing},
  author={Dugast, Florian and Apostolou, Petros and Fernandez, Alfonso and Dong, Wen and Chen, Qian and Strayer, Seth and Wicker, Ryan and To, Albert C},
  journal={Additive Manufacturing},
  volume={37},
  pages={101732},
  year={2021},
  publisher={Elsevier}
}

@article{olleak2022enabling,
  title={Enabling Part-Scale Scanwise process simulation for predicting melt pool variation in LPBF by combining GPU-based Matrix-free FEM and adaptive Remeshing},
  author={Olleak, Alaa and Dugast, Florian and Bharadwaj, Prajwal and Strayer, Seth and Hinnebusch, Shawn and Narra, Sneha and To, Albert C},
  journal={Additive Manufacturing Letters},
  volume={3},
  pages={100051},
  year={2022},
  publisher={Elsevier}
}

@article{liu2023macroscale,
  title={A macroscale thermal simulation strategy with track-scale resolution for laser powder bed fusion},
  author={Liu, Ming and Chiu, Louis NS and Liu, Dedao and Huang, Aijun and Davies, Chris and Wu, Xinhua and Yan, Wenyi},
  journal={Materials \& Design},
  volume={226},
  pages={111622},
  year={2023},
  publisher={Elsevier}
}

@book{strang2007computational,
  title={Computational science and engineering},
  author={Strang, Gilbert},
  volume={791},
  year={2007},
  publisher={Wellesley-Cambridge Press Wellesley, MA}
}

@article{leonor2024go,
  title={GO-MELT: GPU-optimized multilevel execution of LPBF thermal simulations},
  author={Leonor, Joseph P and Wagner, Gregory J},
  journal={Computer Methods in Applied Mechanics and Engineering},
  volume={426},
  pages={116977},
  year={2024},
  publisher={Elsevier}
}

@article{elahi2025efficient,
  title={Efficient part-scale thermal modeling of laser powder bed fusion via a multilevel finite element framework},
  author={Elahi, SM and Leonor, JP and Wu, RY and Wagner, GJ},
  journal={Additive Manufacturing},
  pages={104897},
  year={2025},
  publisher={Elsevier}
}

@article{yang2026efficient,
  title={Efficient thermal simulation in metal additive manufacturing via semi-analytical isogeometric analysis},
  author={Yang, Yang and Ji, Ye and M{\"o}ller, Matthias and Ayas, Can},
  journal={Computer Methods in Applied Mechanics and Engineering},
  volume={457},
  pages={118992},
  year={2026},
  publisher={Elsevier}
}

@article{honarmandi2021rigorous,
  title={A rigorous test and improvement of the Eagar-Tsai model for melt pool characteristics in laser powder bed fusion additive manufacturing},
  author={Honarmandi, P and Seede, R and Xue, L and Shoukr, D and Morcos, P and Zhang, B and Zhang, C and Elwany, A and Karaman, I and Arroyave, R},
  journal={Additive Manufacturing},
  volume={47},
  pages={102300},
  year={2021},
  publisher={Elsevier}
}

@article{zhang2022review,
  title={A review on modelling and simulation of laser additive manufacturing: Heat transfer, microstructure evolutions and mechanical properties},
  author={Zhang, Zhao and Wang, Yifei and Ge, Peng and Wu, Tao},
  journal={Coatings},
  volume={12},
  number={9},
  pages={1277},
  year={2022},
  publisher={MDPI}
}

@article{wolfer2019fast,
  title={Fast solution strategy for transient heat conduction for arbitrary scan paths in additive manufacturing},
  author={Wolfer, Alexander J and Aires, Jeremy and Wheeler, Kevin and Delplanque, Jean-Pierre and Rubenchik, Alexander and Anderson, Andy and Khairallah, Saad},
  journal={Additive Manufacturing},
  volume={30},
  pages={100898},
  year={2019},
  publisher={Elsevier}
}

@article{zagade2021analytical,
  title={Analytical estimation of fusion zone dimensions and cooling rates in part scale laser powder bed fusion},
  author={Zagade, P and Gautham, BP and De, A and DebRoy, T},
  journal={Additive Manufacturing},
  volume={46},
  pages={102222},
  year={2021},
  publisher={Elsevier}
}

@article{ji2022analytical,
  title={Analytical modeling of temperature evolution in laser powder bed fusion considering the size and shape of the build part},
  author={Ji, Xia and Wang, Yongfu and Liang, Steven Y},
  journal={Journal of Materials Processing Technology},
  volume={301},
  pages={117452},
  year={2022},
  publisher={Elsevier}
}

\end{document}